\documentclass[11pt]{article}

\usepackage{graphicx,color}
\usepackage{epsfig}

\usepackage{amsmath,amssymb}
\usepackage{amssymb,eucal,amsmath}

\newcommand{\rank}{\mbox{\rm{~rank~}}}
\newcommand{\diag}{\mbox{diag }}
\newcommand{\Image}{\mbox{\rm{~Image~}}}
\newcommand{\Rww}{\R^{w \times w}}
\newcommand{\Rdotwxi}{\R^{\bullet \times w}[\xi]}

\newcommand{\Camlibel}{\c{C}aml{\i}bel}

\newcommand{\Cinf}{\mathfrak{C}^{\infty}}
\newcommand{\CinfRRw}{\Cinf(\R,\Rw)}

\newcommand{\CinfRRn}{\Cinf(\R,\Rn)}

\newcommand{\Rw}{\R^{\wt}}

\newcommand{\D}{\mathfrak{D}}
\newcommand{\Cn}{\C^{n}}
\newcommand{\Cnn}{\C^{n\times n}}
\newcommand{\Rn}{\R^{n}}

\newcommand{\ze}{\zeta,\eta}
\newcommand{\Rwwze}{\Rww[\ze]}

\newcommand{\Bsub}{{\B_{\rm sub}}}
\newcommand{\Gsub}{{G_{\rm sub}}}

\newcommand{\n}{\tt{n}}

\renewcommand{\S}{\Sigma}

\newcommand{\m}{{\tt{m}}}

\newtheorem{theorem}{Theorem}[section]
\newtheorem{definition}[theorem]{Definition}
\newtheorem{problem}[theorem]{Problem}
\newtheorem{corollary}[theorem]{Corollary}
\newtheorem{prop}[theorem]{Proposition}

\newcommand{\barlambda}{\: \overline{\lambda}}

\newcommand{\jomega}{{j\omega}}
\newcommand{\jj}{{j}}

\newcommand{\EProof}{\hspace*{\fill} {\boldmath $\square$ } \bigskip }

\newcommand{\hinf}{{\mathcal{H}_\infty}}

\newcommand{\Monl} {M(\der)\ell}

\newcommand{\deri}[1]{\mbox{$\frac{{\mathrm d}^{#1}}{{\mathrm d}t^{#1}}$}}
\newcommand{\der}{\deri{}}

\newtheorem{lemma}[theorem]{Lemma}

\newenvironment{proof1}{\noindent
          {\bf Proof.}\ }{\hfill $\Box$ \vspace{0.5cm}\newline}
\newtheorem{exa}[theorem]{Example}

\newcommand{\R}{{\mathbb{R}}}

\newcommand{\C}{{\mathbb{C}}}

\newcommand{\B}{{\mathfrak{B}}}

\newcommand{\Lw}{{\mathfrak{L}^{\tt w}}}
\newcommand{\Lwcont}{{\Lw}_{\rm cont}}

\newenvironment{remark}{\noindent \small
        {\bf Remark:}\ }{\normalsize \vspace{0.5cm}}

\newenvironment{example}{\begin{exa}}{\hfill $\Box$ \end{exa}}

\newcommand{\wt}{\mbox{\rm \texttt{w}}}
\newcommand{\kt}{\texttt{k}}

\newcommand{\Jpq}{J_{q,p}}

\newcommand{\QPsiell}{Q_{\Psi_\ell}}
\newcommand{\QPsix}{Q_{\Psi_x}}
\newcommand{\QPsiw}{Q_{\Psi_w}}
\newcommand{\Bcont}{\B_{\rm cont}}
\newcommand{\Bfullell}{\B_{\rm full}^{\ell}}
\newcommand{\Bfullx}{\B_{\rm full}^{x}}

\usepackage{color}

\title{{\Large \bf Dissipative systems: uncontrollability, observability and RLC
realizability}}
\author{Karikalan Selvaraj, 
Madhu N. Belur and Rihab Abdulrazak\thanks{The authors
are with the Indian Institute of
Technology Bombay, Mumbai 400076, India. Corresponding author email: {\tt
belur@ee.iitb.ac.in}, Fax: +91.22.2572.3707}} %

\begin{document}

\maketitle

\begin{abstract}

The theory of dissipativity has been primarily developed for 
controllable systems/behaviors. 
For various reasons, in the context of uncontrollable systems/behaviors,
a more appropriate definition of dissipativity is in terms of the 
dissipation inequality, namely the {\em existence} of a storage function.
A storage function is a function such that along every
system trajectory, the rate of increase of the storage function
does not exceed the power supplied. While the power supplied is
always expressed in terms of only the external variables, whether or not the 
storage function should
be allowed to depend on only the external variables  and their derivatives
or also unobservable/hidden variables has various consequences on the
notion of dissipativity: this paper
thoroughly investigates the key aspects of both cases, and also proposes
another intuitive definition of dissipativity.

We first assume that the storage function can be expressed in terms
of the external variables and their derivatives only and prove our
first main result that, assuming the uncontrollable poles
are unmixed, i.e. 
no pair of uncontrollable poles add to zero, and assuming a
strictness of dissipativity at the infinity frequency, the dissipativities
of a system and its controllable part are equivalent; in other words
once the autonomous subsystem satisfies a Lyapunov equation solvability-like
condition, it does not interfere with the dissipativity of the system.
We also show that the storage function in this case is a static
state function. 
This main result proof involves new results about solvability of the
Algebraic Riccati Equation, and uses techniques from Indefinite
Linear Algebra and Hamiltonian matrix properties.

We then investigate the utility of unobservable/hidden variables
in the definition of storage function: we prove that
lossless uncontrollable behaviors are ones which require
storage function to be expressed in terms of variables that are
unobservable from the external variables.

We next propose another intuitive definition: a behavior is called
dissipative if it can be embedded in a 
controllable dissipative {\em super-behavior}. We show that this 
definition imposes a constraint on the number of inputs and thus
explains unintuitive examples from the literature in the context
of lossless/orthogonal behaviors.  These results are finally 
related to RLC realizability of passive networks, specifically
to the nonrealizability of the nullator one-port circuit using RLC elements.
\end{abstract}

\section{Introduction}
The theory of dissipativity for linear dynamical systems helps in the analysis 
and design of control systems for several control problems,  for example,
LQR/LQG control, $\mathcal{H}_{\infty}$, synthesis of 
passive systems, and optimal estimation problems. When dealing
with LTI systems, it is straightforward
to define dissipativity for controllable systems due to a certain property
of such systems that their compactly supported system trajectories are,
loosely speaking, `dense' in the set of all allowed trajectories.
However, this is not the case for uncontrollable systems, and this
situation is the central focus of this paper.  We elaborate more on
this point when we define dissipativity and review equivalent conditions
for controllable systems in Section \ref{sec2}.

In this paper, we use a less-often-used definition of dissipativity
for systems, possibly uncontrollable, and generalize key results
using some techniques from indefinite linear algebra (see \cite{GLR05}) for
solving Algebraic Riccati Inequalities in the context of an uncontrollable
state space system. Like in \cite{CWB03,PB08}, we define a system
as dissipative if there exists a storage function that satisfies
the dissipation inequality for all system trajectories. The existential
aspect of this definition raises key issues that this paper deals with.

The main result we show is that if the uncontrollable poles of an
LTI system are such that no two of them add to zero, and if the
controllable subsystem strictly dissipates energy at frequency equal
to infinity, then the
dissipativity of the controllable subsystem is equivalent to the
system's dissipativity. 
We also show that, using the concatenability
axiom of the state, the energy stored in a system is a static function
of the state variables. Further, we also show that
this state is `observable' from the external variables, i.e.
the state is a linear combination of the external variables and possibly
their derivatives. This is intuitively expected in view of the fact that energy
exchange between the system and its ambience
takes place through the external variables. However, it appears
that this may not be the case for lossless systems, i.e. systems
that don't dissipate any energy, nor contain a source within.

\section{Preliminaries} \label{sec2}

In this section we include various definitions about the behavioral
framework for studying dynamical systems (Subsection \ref{subsec2.1}) and
then introduce background results about dissipative
systems (Subsection \ref{subsec2.2}). Subsection \ref{subsec2.3} contains
brief notation about indefinite linear algebra from \cite{GLR05}. For this
paper, $\R$ denotes the set of all real numbers and $\R[\xi]$ the set
of polynomials in the indeterminate $\xi$ and real coefficients; matrices
and polynomial matrices are denoted the obvious way. We use $\bullet$ to
leave a row dimension unspecified, for example, $\Rdotwxi$. The
space of infinitely often differentiable functions from $\R$ to say
$\Rn$ is denoted by $\CinfRRn$ and $\D$ denotes the set of compactly supported
functions within this space.

\subsection{The behavioral approach} \label{subsec2.1}

When dealing with linear differential systems, it is convenient to
use polynomial matrices for describing a differential equation.
Suppose $R_0$, $R_1 \dots$ $R_N$ are constant matrices of the same size such that
\[
R_0 w + R_1 \der w + R_2 \deri{2} w \dots R_N \deri{N} w =0
\]
is a linear constant coefficient ordinary differential equation in
the variable $w$.
We define the polynomial matrix $R(\xi):= R_0 + R_1 \xi + \dots
R_N \xi^N $, and represent the above differential equation as $R(\der)w=0$.

A linear differential behavior, denoted by $\mathfrak{B}$, is defined
as the set of 
all infinitely often differentiable trajectories
that satisfy a system of ordinary 
linear differential equations with constant coefficients, i.e.,
\begin{equation*}
\B := \{w \in \mathfrak{C}^\infty(\mathbb{R},\mathbb{R}^\texttt{w})
\: \mid \: R(\frac{d}{dt})w = 0 \},
\end{equation*}
where $R(\xi)$ is a polynomial matrix having $\texttt{w}$ number of columns, 
i.e., $R(\xi)\in \Rdotwxi$. We 
denote the set of all such linear differential behaviors with 
$\texttt{w}$ number of variables by $\Lw$. The linear 
differential behavior $\mathfrak{B} \in \mathfrak{L}^\texttt{w}$ given by the 
above definition can also be represented by $\B = 
\ker{R(\frac{d}{dt})} $. One views the rows of $R$ as differential equations
that the variable $w$ has to satisfy in order for a trajectory $w(t)$ to
be in the behavior $\B$. The matrix $R$ is not unique and one can use elementary
row operations to modify $R$ and this does not change the set
of solutions $\B$: this thus allows assuming without loss of generality
that $R$ has full row rank (see \cite{Wil07-2}). 
The number of inputs of $\B$ is defined as
$\wt-\rank(R)$ and is called the input
cardinality. This integer depends only on $\B$ and not on the $R$ used to
define it; $\texttt{m}(\mathfrak{B})$ denotes the number of inputs.

\newcommand{\mB}{{\texttt{m}(\mathfrak{B})}}

An important fundamental concept is controllability of a system.
A behavior $\mathfrak{B}= \ker{R(\frac{d}{dt})w}$ is said to be 
\emph{controllable}, if for every $w_1$ and $w_2 \in \mathfrak{B}$,
there exist $w_3 \in \mathfrak{B}$ and $\tau >0$ such that
\[
w_3 =
\left\{ \begin{array}{l} w_1(t) \: \: \mbox{for all $t\leqslant 0$}, \\
w_2(t) \: \: \mbox{for all t $\geqslant \tau$}.
\end{array} \right. 
\]
The set of all {\em controllable}
behaviors with $\texttt{w}$ variables is denoted 
as $\Lwcont$. This patchability definition of controllability
is known to have the
traditional Kalman state-space definition of controllability as a special
case in \cite{Wil07-2}. Further, it is shown that $\mathfrak{B} = 
\ker{R(\frac{d}{dt})w} $ is controllable if and only if $R(\lambda)$ has
constant rank for all $\lambda \in \mathbb{C}$.  It is also 
shown in \cite{Wil07-2} that $\mathfrak{B}$ is
controllable if and only if it can be defined as
\begin{equation*}
\begin{array}{c}
\mathfrak{B} := \{w \in \mathfrak{C}^\infty(\mathbb{R},\mathbb{R}^\texttt{w}) 
\: \mid \: \mbox{there exists an} \: \ell \in 
\mathfrak{C}^\infty(\mathbb{R},\mathbb{R}^\texttt{m}) \: \qquad \\
 \qquad \text{such that} \: w = M(\frac{d}{dt})\ell \},
\end{array}
\end{equation*}
where $M(\xi) \in \mathbb{R}^{\texttt{w} \times \texttt{m}}\left[ \xi 
\right]$. This representation of $\mathfrak{B} \in  
\mathfrak{L}^{\texttt{w}}_{cont}$ is known as an \emph{image representation}.
The variable $\ell$ is called a \emph{latent} variable: these are auxiliary
variables used to describe the behavior; we distinguish the variable $w$ as
 the manifest variable, the variable of interest.
It is known that $\mathfrak{B} \in  \mathfrak{L}^{\texttt{w}}_{cont}$ 
always allows an image representation with $M(\xi)$ such that $M(\lambda)$
has full column rank for every $\lambda \in \C$. This 
kind of image representation is known as an \emph{observable} image 
representation. In this paper, unless otherwise stated
explicitly, we assume the image representations are observable.
The use of the term `observable' is motivated by the fact that
the variable $\ell$ is \emph{observable} from the variable $w$.
This notion is defined as follows. 

For a behavior $\B$ with variables $w$ and $\ell$, we say $\ell$ is observable
from $w$ if whenever $(w,\ell_1)$ and $(w,\ell_2)$ both  are in $\B$, we
have $\ell_1=\ell_2$. Observability of $\ell$ from $w$ in a behavior $\B$
implies that there exists a polynomial matrix $F(\xi)$ such
that $\ell=F(\der) w$ for all $w$ and $\ell$ in the behavior.

We now define relevant notions in the context of uncontrollable behaviors.
For a behavior $\mathfrak{B}$, possibly uncontrollable, 
the largest controllable behavior contained in 
$\mathfrak{B}$ is called the \emph{controllable} part of 
$\mathfrak{B}$, and denoted by $\mathfrak{B}_{cont}$.
An important fact about the 
controllable part of $\mathfrak{B}$ is that $\texttt{m}(\mathfrak{B}_{cont}) = 
\texttt{m}(\mathfrak{B})$. 
The set of complex numbers $\lambda$ for which 
$R(\lambda)$ loses rank is called the set of \emph{uncontrollable modes} 
and is denoted by $\Lambda_{un}$. For a detailed exposition on behaviors, 
controllability and observability we refer the reader to \cite{Wil07-2}.

\subsection{Quadratic Differential Forms and dissipativity} \label{subsec2.2}

The concept of Quadratic Differential Forms (QDF) (see \cite{WT98}) 
is central to this paper. 
Consider a two 
variable polynomial matrix with real coefficients, $\Phi(\zeta,\eta) := 
\sum_{j,k}{\Phi_{jk}\zeta^{j} \eta^{k}} \in \mathbb{R}^{\texttt{w}\times 
\texttt{w}}\left[\zeta,\eta \right]$, where $\Phi_{jk} \in 
\mathbb{R}^{\texttt{w}\times \texttt{w}}$,
The QDF $Q_\Phi$ induced by $\Phi$
is a map $Q_\Phi: 
\mathfrak{C}^\infty(\mathbb{R},\mathbb{R}^\texttt{w}) \longrightarrow 
\mathfrak{C}^\infty(\mathbb{R},\mathbb{R})$ defined as 
\begin{equation*}
Q_\Phi(w) := \sum_{j,k} (\frac{d^jw}{dt^j})^T \Phi_{jk} (\frac{d^kw}{dt^k}).
\end{equation*}
When dealing with quadratic forms in $w$ and its derivatives, we can assume 
without loss of generality that $\Phi(\zeta,\eta) = \Phi^T(\eta,\zeta)$:
such a $\Phi$ is called
a symmetric two variable polynomial matrix. 
A quadratic form induced by a 
real symmetric constant matrix $S \in \mathbb{R}^{\texttt{w}\times 
\texttt{w}}$ is a special QDF. We frequently need the number of
positive and negative eigenvalues of a nonsingular, symmetric matrix
$S$: they are denoted by $\sigma_+(S)$ and $\sigma_-(S)$ respectively.

For a two variable polynomial matrix 
$\Phi(\zeta,\eta) \in \mathbb{R}^{\texttt{w}\times \texttt{w}}\left[\zeta,\eta 
\right]$, we define the single variable polynomial matrix $\partial \Phi$
by $\partial \Phi(\xi)=\Phi(-\xi ,\xi)$.

Consider $\Sigma \in \Rww$, a symmetric nonsingular matrix. 
A behavior $\B\in\Lw$ is said to be dissipative with respect to the
supply rate $\Sigma$ (or $\Sigma$-dissipative)
if there exists a QDF $Q_\Psi$ such that
\begin{equation} \label{dissipineq}
\der Q_\Psi(w) \leqslant w^T \S w \mbox{ for all } w \in \B~.
\end{equation}

The QDF $Q_\Psi$ is called a storage function: it signifies the energy
stored in the system at any time instant. The above inequality is
called the dissipation inequality. A behavior $\B$ is called $\Sigma$-lossless
if the above inequality is satisfied with an equality for some QDF $Q_\Psi$.
Notice that the storage function plays the same role as that
of Lyapunov functions in the context of autonomous systems; the notion
of storage functions
is a generalization to non-autonomous systems of Lyapunov functions, as pointed
in \cite{WT98}.  The following theorem from \cite{WT98} applies to controllable
behaviors.
\begin{prop} Consider $\B \in \Lwcont$ and let $w=\Monl$ be an
observable image representation. Suppose $\Sigma\in\Rww$ is symmetric
and nonsingular. Then, the following are equivalent.
\begin{enumerate}
\item There exists a QDF $Q_\Psi$ such that inequality \eqref{dissipineq}
 is satisfied for all $w\in\B$.
\item $\int_\R w^T \S w dt \geqslant 0$ for all $w\in \B\cap \D$, the
compactly supported trajectories in $\B$.
\item $M^T(-\jomega) \Sigma M(\jomega) \geqslant 0$ for all $\omega \in \R$.
\end{enumerate}
\end{prop}

The significance of the above theorem is that, for controllable
systems, it is possible to verify dissipativity by checking
non-negativity of the above integral over all compactly supported
trajectories: the compact support signifying that we calculate the
`net power' transferred when the system starts `from rest' and `ends at rest'.
The starting and ending `at rest' ensures that for linear
systems there is no internal energy at this time. The absence of internal
energy allows ruling out the storage function from this condition: in fact, this
is used as the definition of dissipativity for controllable systems.
The same cannot be done for uncontrollable systems due to the
compactly supported trajectories not being dense in the behavior (see 
\cite{PilSha:98}).
An extreme case is an autonomous behavior, i.e. a behavior which has
$\mB=0$: while the zero trajectory
is the only compactly supported trajectory,
the behavior consists of exponentials corresponding to the uncontrollable
poles of $\B$.  The issue of existence of storage functions
is elaborated in \cite[Remark 5.9]{WT98} and
in text following \cite[Proposition 3.3]{PB08}.

\subsection{States} \label{subsec2.3}

The state variable is defined as a latent variable 
that satisfies the \emph{property of state}, that is, 
if $(w_1,x_1)$, $(w_2,x_2) \in \mathfrak{B}_f$ and $x_1(0)=x_2(0)$,
then the new 
trajectory $(w,x)$ formed by
concatenating $(w_1,x_1)$ and $(w_2,x_2)$ at $t=0$, 
i.e., 
$$
(w,x)(t) =
\begin{array}{l}
 (w_1,x_1)(t) \: \: \: \mbox{ for all } \quad t\leqslant 0 \\
 (w_2,x_2)(t) \: \: \: \mbox{ for all } \quad  t > 0,
\end{array}
$$
also satisfies the system equation of $\mathfrak{B}_f$ in a distributional 
sense \cite{WT02}. It is intuitively expected that a variable $x$ has the state
property if and only if $w$ and $x$ satisfy an equation that is
at most first order in $x$ and zeroth order in $w$:
see \cite{RW97} for precise statement formulation and proof. When
$w$ is partitioned into $w=(w_1,w_2)$, with $w_1$ as the input and $w_2$ as
the output, then $\mathfrak{B}_f$ admits the more familiar
\emph{input/state/output} (i/s/o) representation as 
\begin{equation} \label{eq:state}
\frac{d}{dt}x = Ax + Bw_1, \: \: \: \: w_2 = Cx + Dw_1.
\end{equation}
One can ensure that $x$ is observable from $w$; this is equivalent
to conventional observability of the pair $(C,A)$. While such a state space
representation is admitted by any $\B$,
the pair $(A,B)$ is controllable (in the state space sense) if and only if 
$\mathfrak{B}$ is controllable (in the behavioral sense defined above).

\subsection{Indefinite linear algebra}

In this paper, we use certain properties of matrices that are self-adjoint
with respect to an {\em indefinite} inner product. 
We briefly review self-adjoint matrices and neutral 
subspaces (see \cite{GLR05}). Let $P \in \Cnn$ be an invertible 
Hermitian matrix. This defines an indefinite inner product on $\Cn$ by
$(Px,x):=x^* P x$, where $x^*$ is the complex conjugate transpose of the
vector $x\in\Cn$. For a complex matrix $A$, the complex conjugate transpose
of $A$ is denoted by $A^*$.

Consider matrices $A$ and $P\in\Cnn$ with $P$ invertible and Hermitian.
The $P$-adjoint of the matrix $A$, denoted by $A^{\left[\ast \right]}$,
is defined as 
$A^{\left[\ast \right]} = P^{-1}A ^* P$. The matrix $A$ is
said 
to be $P$-self-adjoint 
if $A = A^{[\ast]}$, i.e., $A = P^{-1}A ^* P$. 

Since $P$ is not sign-definite in general, the sign of $(Px,x)$ is zero,
positive or negative depending on the vector $x$.
A subspace $\mathcal{M} \subseteq \mathbb{C}^n$ is said to be $P$-neutral if 
the inner product $(Px,x) = 0$ for all $x \in \mathcal{M}$.

\section{Dissipativity of uncontrollable behaviors}

In this paper, we deal with systems which satisfy the dissipativity
property, i.e. \emph{net} energy is directed \emph{inwards}
along every system trajectory. As elaborated below, the `total' aspect
of the energy involves an integral, thus bringing in the initial and
final conditions of the trajectory being integrated. The 
convenience of starting-from-rest and ending-at-rest applies
to only controllable systems as we will review soon. The notion of
storage function helps in formulating the dissipation property as
an inequality to be satisfied {\em at each time-instant}.
A central issue in this paper is what variables
should the storage function be allowed to depend on. We explore
dependencies on a latent variable $\ell$,
or on a state variable $x$, or on the manifest variable $w$: this 
is indicated in the storage function subscript; the following definition
has appeared in several works, see \cite{WT98, Wil04, Wil07-2}, for example.

\begin{definition} \label{defn:dissip}
Let $\Sigma \in \Rww$ be a nonsingular symmetric matrix, inducing
the supply rate $w^T \Sigma w$.
Consider a behavior $\B \in \Lw$ with manifest variables
$w$ and latent variable $\ell$, with the
corresponding full behavior $\Bfullell$.
For the behavior $\B$, let $x$ be a state variable with
the corresponding full behavior $\Bfullx$. 
With respect to the supply rate $\Sigma$, the behavior $\B$ is said to be 
dissipative if there exists a quadratic differential form $Q_\Psi$ such that
\begin{equation} \label{eq:dissip}
\der \QPsiell(\ell) \leqslant w^T \Sigma w
\mbox{ for all $(w,\ell) \in \Bfullell$.}
\end{equation} 
\begin{enumerate}
\item 
The function $\QPsiell$, a quadratic function of $\ell$ and its
derivatives, is called a storage function.
\item A storage function $\QPsiell$ is said
to be an {\em observable} storage function if the latent variable
$\ell$ is observable from the manifest variable $w$. In this case,
there exists a storage function $\QPsiw$ such
that $\der \QPsiw(w)\leqslant w^T \Sigma w$ for all $w\in \B$.
\item A storage function $\QPsix(x)$ is said to be a \emph{state function}
if $\QPsix(x)$ is equal to $x^T K x$ for some constant matrix $K$.
\end{enumerate} 
\end{definition}

In this paper, we study dissipativity with respect 
to a constant, nonsingular, symmetric 
$\Sigma \in \mathbb{R}^{\texttt{w} \times \texttt{w}}$. It is known (see
\cite[Remark 5.11]{WT98} and \cite[Proposition 2]{WT02}) that
$\Sigma$-dissipativity of a behavior $\B$ implies that the input cardinality
of $\B$ cannot exceed $\sigma_+(\Sigma)$, i.e.,
$\texttt{m}(\mathfrak{B}) \leqslant \sigma_+ (\Sigma) $.  In this context
it is helpful to perform a coordinate transformation in the variables
$w$ so that $\Sigma$ is a diagonal matrix consisting of only $+1$'s and
$-1$'s along the diagonal. Moreover, there exists an input/output
partition such that all the inputs correspond to $+1$'s only  and such that
the transfer function from these inputs to all other variables is {\em proper}
(see \cite[Remark 5.11]{WT98}). In view of these facts and the inequality
$\texttt{m}(\mathfrak{B}) \leqslant \sigma_+ (\Sigma) $, we assume
without loss of generality 
\begin{equation} \label{eq:Sigmanew}
\Sigma = 
\begin{bmatrix}
	I_{\texttt{m}} & 0 & 0 \\
	0 & I_{\texttt{q}} & 0 \\
	0 & 0 & -I_{\texttt{p}} \\
\end{bmatrix} \: \: \: \: \text{and define} \: \: \: \: 
J_{pq} = \left[ \begin{array}{ccc}
	I_{\texttt{q}} & 0 \\
	0 & -I_{\texttt{p}} \\
\end{array} \right]
\end{equation}
where $\texttt{m}$ is the number of inputs in $\B$.

\subsection{Dissipativity of uncontrollable behaviors: main results}

The following theorem assumes an unmixing condition - no pair of
 uncontrollable poles are symmetric with respect to the imaginary axis. If
 this unmixing condition is satisfied for the uncontrollable poles, then the
 controllable part of a behavior being dissipative is equivalent to the
 dissipativity of the whole behavior. 

\begin{theorem}
Consider a behavior $\mathfrak{B} \in \mathfrak{L}^\texttt{w}$ and a nonsingular,
 symmetric $\Sigma \in \mathbb{R}^{\texttt{w} \times \texttt{w}}$	with the
 input cardinality of $\mathfrak{B}$ at-most the positive signature of
 $\Sigma$, i.e., $\texttt{m}(\mathfrak{B}) \leqslant \sigma_+ (\Sigma) $.
 Assume that the uncontrollable poles are such that 
$\Lambda_{un} \cap -\Lambda_{un} = \emptyset $.
Let $\mathfrak{B}$ have an observable image representation, 
$\mathfrak{B} = \Image M(\frac{d}{dt})$, where $M(\xi)$ is partitioned as 
\begin{equation}
M(\xi) = 
\begin{bmatrix} W_1(\xi) \\ W_2(\xi) \end{bmatrix}; \: 
\: \: \: W_1 \in \mathbb{R}^{\texttt{m} \times \texttt{m}}\left[ \xi \right],
 W_2 \in \mathbb{R}^{(\texttt{p+q}) \times \texttt{m}}\left[ \xi \right].
\end{equation}
Let $G(\xi) := W_2(\xi)W_1(\xi)^{-1}$ and $D := \lim_{\omega \to \infty}
 G(j\omega)$. Assume $\mathfrak{B}_{cont}$ is such that $(I_m+D^TJ_{pq}D) > 0$.
Then, $\mathfrak{B}$ is $\Sigma$-dissipative if and only if its controllable
 part $\mathfrak{B}_{cont}$ is $\Sigma$-dissipative.
\label{main}
\end{theorem}

The above result gives conditions under which the autonomous part of
a behavior plays no hindrance to dissipativity of the behavior after the
controllable part is dissipative. One of the conditions for this
is that the uncontrollable poles are not `mixed', meaning no two of
the uncontrollable poles add to zero. For autonomous LTI systems, this
condition is a necessary and sufficient condition for solvability of
the Lyapunov equation. Of course, storage functions are just
generalizations of Lyapunov functions to non-autonomous systems.
The other condition: $(I_m+D^TJ_{pq}D) > 0$ on the controllable part
is a kind of strictness of dissipativity `at the
infinity\footnote{The feed-through term $D$ is finite, since the transfer
function is proper; see text before Equation \eqref{eq:Sigmanew}.} frequency'.
While this condition allows the use of Hamiltonian matrices in the proofs,
this condition also rules out consideration of
lossless systems from the above result. More significance of
the assumption is noted in the remark below.

\begin{remark}
While a behavior $\B\in\Lw$ admits
many i/o partitions, and for each i/o partition,
admits many i/s/o representations, the matrix $I+D^T \Jpq D$ depends
only on the behavior, and in particular, only on $\Bcont$ the controllable
part. In other words, this matrix can be found from any image representation
of $\Bcont$. This matrix being positive
semi-definite is a necessary condition for dissipativity of $\Bcont$, and
this denotes dissipativity at very high frequencies, i.e. as
$\omega \rightarrow \infty$. Positive definiteness of $I+D^T \Jpq D$
may hence be termed as `strict dissipativity at the infinity frequency'.
This assumption helps in the existence of a Hamiltonian matrix, our proofs
use the Hamiltonian matrix properties and techniques from indefinite
linear algebra (see \cite{GLR05}).
Positive definiteness of
this matrix is guaranteed, for example, by strict dissipativity of
a behavior; on the other hand, lossless controllable behaviors
have this matrix as zero. Another interpretation of 
the condition $I+D^T \Jpq D > 0$ is that the `memoryless' part of the
behavior is strictly dissipative; the memoryless/static part of
a behavior was defined in \cite{WT02}, and we don't require this notion
in this paper.
\end{remark}

We now review some existing results and formulate/prove
new results about Hamiltonian matrices
in the context of dissipativity of controllable and uncontrollable systems. 
These are required for the proof of the main result.

To prove that the behavior $\mathfrak{B}$ is $\Sigma$-dissipative, we use the 
proposition below and show that the existence of a symmetric solution to the 
dissipation LMI is sufficient for $\Sigma$-dissipativity. Let the behavior 
have an input/output partition $w=(w_1,w_2)$, $w_1 \in 
\mathfrak{C}^{\infty}(\mathbb{R},\mathbb{R}^\texttt{m})$ as the input and $w_2 
\in \mathfrak{C}^{\infty}(\mathbb{R},\mathbb{R}^\texttt{p+q})$ as the output. 
Let the input/state/output representation of the behavior $\B$ be,
\begin{equation} \label{eq:iso3}
\frac{d}{dt}x = Ax + Bw_1, \: \: \: \: w_2 = Cx + Dw_1.
\end{equation}
with $(C,A)$ observable. The following result is well-known: see
\cite{Wil:72, Sch:90, BoyGhaFerBal:94}.

\begin{prop}\label{LMIunctrb}
Let a behavior $\mathfrak{B} \in \mathfrak{L}^{\texttt{w}}$ have an i/s/o 
representation for an input/output partition $w=(w_1,w_2)$ with $A,B,C,D$ as 
state space matrices and $(A,B)$ possibly uncontrollable. The behavior 
$\mathfrak{B}$ is $\Sigma$-dissipative if there exists a symmetric $K \in 
\mathbb{R}^{\texttt{n}\times \texttt{n}}$ which solves the following LMI
\begin{equation} \label{eq:LMInew} 
\left[
\begin{array}{cc}
	(KA + A^TK - C^TJ_{pq}C) & (KB - C^TJ_{pq}D) \\
	(KB - C^TJ_{pq}D)^T & -(I_\texttt{m} + D^TJ_{pq}D) \\
\end{array}
\right]
\leqslant 0.
\end{equation}
Further, under the condition that $(A,B)$ is controllable, existence
of $K$ solving the above LMI is necessary and sufficient for $\S$-disspativity
of $\B$.  
\end{prop}


As $(I_\texttt{m}+D^TJ_{pq}D)$ is invertible, the Schur complement of 
$(I_\texttt{m}+D^TJ_{pq}D)$ in the above LMI gives the Algebraic Riccati 
Inequality (ARI)
\begin{equation} \label{eq:ARInew}
\begin{array}{c}
K(A - B(I_\texttt{m}+D^TJ_{pq}D)^{-1}D^TJ_{pq}C)\\
 + (A - B(I_\texttt{m}+D^TJ_{pq}D)^{-1}D^TJ_{pq}C)^TK  \\ 
+ KB(I_\texttt{m}+D^TJ_{pq}D)^{-1}B^TK - C^T(J_{pq}+DD^T)^{-1}CA 
\end{array}
 \: \: \: \: \leqslant 0
\end{equation}
The corresponding equation is the Algebraic Riccati Equation (ARE) and we use 
properties of this ARE in proving Theorem \ref{main}. We rewrite the
above ARE as 
\begin{equation} \label{eq:ARE} 
K\tilde{A} +\tilde{A}^TK + K\tilde{D}K - \tilde{C} = 0
\end{equation}
where
\begin{equation} \label{eq:tildenew}
\tilde{A}:= (A - B(I_\texttt{m}+D^TJ_{pq}D)^{-1}D^TJ_{pq}C), 
\nonumber
\end{equation}
and
\begin{equation*}
\tilde{D}:= B(I_\texttt{m}+D^TJ_{pq}D)^{-1}B^T, ~ ~  \tilde{C}:= 
C^T(J_{pq}+DD^T)^{-1}C.
\end{equation*}

Also define the Hamiltonian matrix corresponding to the ARE: 
\begin{equation} \label{eq:ham}
	H = \left[ \begin{array}{cc}
				\tilde{A} & \tilde{D} \\
				\tilde{C} & -\tilde{A}^\ast \\ 
			\end{array} \right].
\end{equation}

We define the following matrices; they play a crucial role in the results
we use from \cite{GLR05} and in our proofs.
\begin{equation} \label{eq:ham1}
	M := \jj \:H,  \: \: \: \: 
	P := \left[ \begin{array}{cc}
				-\tilde{C} & \tilde{A}^\ast \\
				\tilde{A} &  \tilde{D}\\ 
			\end{array} \right] \: \: \: \text{and} \: \: \:
	\widehat{P} := \jj \: \left[ \begin{array}{cc}
				0 & I \\
				-I &  0\\ 
			\end{array} \right].
\end{equation}


It is well-known that a symmetric solution to ARE can be obtained from an 
$n$-dimensional, $M$-invariant, $P$-neutral subspace and we state this in the 
proposition below. Let $K \in \mathbb{C}^{\n\times \n}$ be a Hermitian 
matrix. The graph subspace corresponding to matrix $K$ is defined as,
\begin{equation*}
\mathcal{G}(K) := \text{Im} \left[ \begin{array}{c}
	I \\
	K \\
\end{array}  \right]. 
\end{equation*}
Notice that $\mathcal{G}(K)$ is an n-dimensional subspace of $\mathbb{C}^{2n}$.
The following proposition from \cite{GLR05} states solvability of the ARE
in terms of $P$-neutrality with respect to the above $P$.

\begin{prop}\label{graphneutral}(See \cite{GLR05})
Consider the ARE \eqref{eq:ARE} with $(\tilde{A},\tilde{D})$ possibly 
uncontrollable. Then $K \in \mathbb{C}^{\n \times \n}$ is a Hermitian 
solution of the ARE if and only if the graph subspace of $K$ is $M$-invariant 
and $P$-neutral. 
\label{pro:invsln}
\end{prop}

The proposition below allows use of simplified $A$ and $B$ for all
later purposes when dealing with uncontrollable systems.

\begin{prop}(See \cite{Won79})
Consider the behavior $\mathfrak{B} \in \mathfrak{L}^{\texttt{w}}$ with an 
input/state/output representation
$\frac{d}{dt}x = Ax + Bw_1$, $w_2 = Cx + Dw_1$, where $w = (w_1,w_2)$. Then 
there exists a nonsingular matrix $T \in \mathbb{R}^{\n \times \n}$ such that 
\begin{equation*}
T^{-1}AT = \left[ \begin{array}{cc}
A_c & A_{cp} \\
0 & A_u \\
\end{array} \right], \: \: \: T^{-1}B = \left[ \begin{array}{c}
B_c \\
0 \\
\end{array} \right] \: \: and \: \: CT = \left[ \begin{array}{cc}
C_c & C_u \\
\end{array} \right].
\end{equation*}
Further,
\begin{equation*}
\frac{d}{dt}x = A_cx + B_cw_1,\: \: \: \:  w_2 = C_cx + Dw_1.
\end{equation*}
gives an i/s/o representation for the controllable part 
$\mathfrak{B}_{cont}$.
\end{prop}

Let $H_c$ and $M_c$ denote the corresponding matrices for the controllable 
part $\mathfrak{B}_{cont}$ as defined in \eqref{eq:ham1} with $\tilde{A}_c, 
\tilde{D}_c$ and $\tilde{C}_c$ defined accordingly in \eqref{eq:tildenew}.
	
\begin{lemma}\label{partial}
Suppose $\mathfrak{B}_{cont} \in \mathfrak{L}^{\texttt{w}}_{cont}$ satisfies 
the assumption that $(I_m+D^TJ_{pq}D) > 0$. If $\mathfrak{B}_{cont}$ is 
$\Sigma$-dissipative, then the partial multiplicities corresponding to the 
real eigenvalues of $M_c$, if any, are all even.
\end{lemma}

In order to prove Lemma \ref{partial}, we use a
result from \cite{GLR05} concerning the 
partial multiplicities of real eigenvalues of $M$. Let the generalized 
eigenspace of a matrix $A$ corresponding to an eigenvalue $\lambda_0$ be 
denoted by $\mathcal{R}_{\lambda_0}(A)$ and the controllable subspace of the 
pair $(\tilde{A},\tilde{D})$ be denoted by $\mathfrak{C}_{\tilde{A},\tilde{D}}$.

\begin{prop}\label{partialM}(See \cite{GLR05})
Consider the behavior $\mathfrak{B}_{cont} \in 
\mathfrak{L}^{\texttt{w}}_{cont}$. Assume $\tilde{D}_c \geq 0$ and 
$\tilde{C}_c^{\ast} = \tilde{C}_c$ and there exists a Hermitian solution $K 
\in \mathbb{C}^{\n \times \n}$ to the ARE \eqref{eq:ARE}.
Suppose 
\begin{equation*}
\mathcal{R}_{\lambda_0}(\tilde{A}_c+\tilde{D}_cK) \subseteq 
\mathfrak{C}_{\tilde{A}_c,\tilde{D}_c}
\end{equation*}
for every purely imaginary eigenvalue $\lambda_0$ of 
$(\tilde{A}_c+\tilde{D}_cK)$. Then, the partial multiplicities of corresponding 
real eigenvalues of $M_c$ are all even and are twice the partial 
multiplicities of those corresponding to the purely imaginary eigenvalues of 
$(\tilde{A}_c+\tilde{D}_cK)$.
\end{prop}

\noindent 
\textbf{Proof of the Lemma \ref{partial}}: As the controllable part is 
$\Sigma$-dissipative, there exists a symmetric solution $K$ to the ARE by 
Proposition \ref{LMIunctrb}. Since the behavior is controllable, the following 
is true: $\mathfrak{C}_{\tilde{A}_c,\tilde{D}_c} = \mathbb{C}^n$ and hence 
$\mathcal{R}_{\lambda_0}(\tilde{A}_c+\tilde{D}_cK) \subseteq 
\mathfrak{C}_{\tilde{A}_c,\tilde{D}_c}$ for every purely imaginary eigenvalue 
$\lambda_0$ of $(\tilde{A}_c+\tilde{D}_cK)$. Thus, using
Proposition \ref{partialM}, the partial multiplicities of all real
eigenvalues of $M_c$, 
if any, are all even. This completes the proof.

We define a set called \emph{c-set} as in \cite{GLR05}. Such a \emph{c-set}, 
if exists, guarantees the existence of a unique $P$-neutral, $M$-invariant 
subspace $\mathcal{N}$ under certain conditions. This is made precise in the 
proposition below.

\begin{definition}(See \cite{GLR05})
Let $M\in\Cnn$ and
let $\mathcal{C}$ be a finite set of non-real complex numbers. $\mathcal{C}$ 
is called a \underline{c-set} of $M$ if it satisfies the following properties
\begin{enumerate}
	\item $\mathcal{C} \cap \overline{\mathcal{C}} = \emptyset$ 
	\item $\mathcal{C} \cup \overline{\mathcal{C}}$ = $\sigma(M) \setminus 
\mathbb{R}$, the set of all non-real eigenvalues of $M$.
\end{enumerate}
\end{definition}

\newcommand{\CsetC}{{\mathcal{C}}}

Let $\mathcal{N}$ be an invariant subspace of $M$ and denote the restriction 
of $A$ to $\mathcal{N}$ as $A|_{\mathcal{N}}$.

\begin{prop}\label{invariant}(See \cite{GLR05})
Let $M \in \mathbb{C}^{2n \times 2n}$ be a $P$-self-adjoint matrix such that 
the sizes of the Jordan blocks of $M$, say: $m_1,m_2,...,m_r$, corresponding
to real 
eigenvalues of $M$ are all even.
Then for every c-set $\CsetC$ there exists a unique 
$P$-neutral $M$-invariant subspace $\mathcal{N}$ of dimension $n$ and 
$\sigma(M|_{\mathcal{N}}) \setminus \mathbb{R} = \mathcal{C}$, and the sizes 
of the Jordan blocks of $M|_{\mathcal{N}}$ corresponding to the real 
eigenvalues are $\frac{1}{2}\,m_1,\frac{1}{2}\,m_2,...,\frac{1}{2}\,m_r$.
\end{prop}

It is easy to verify that $M$ and $P$ satisfy the following relation, $PM = 
M^\ast P$. Hence $M$ is $P$-self-adjoint and Proposition \ref{invariant} can be 
used.


Following proposition, which is a reformulation and combination
of Theorems A.6.1, A.6.2 and A.6.3 in \cite{GLR05},
states that the partial multiplicities of an eigenvalue 
are unaffected by pre-multiplying and/or post-multiplying by an unimodular 
matrix. This result is used in the proof of the Theorem \ref{main}.

\begin{prop}\label{evenparmul}
Consider $S_1(\xi)$ and $S_2(\xi) \in \mathbb{R}^{\texttt{w} \times 
\texttt{w}} \left[ \xi\right]$ and let $p_j^1(\xi)$ and $p_j^2(\xi) \in 
\mathbb{R}\left[ \xi\right]$ for $i=1,\dots,w$ be the invariant polynomials of 
$S_1(\xi)$ and $S_2(\xi)$ respectively. Suppose $S_1(\xi) = T_1 S_2(\xi) T_2$, 
for invertible $T_1, T_2 \in \mathbb{R}^{\texttt{w} \times \texttt{w}}$. Let 
$\lambda \in \mathbb{C}$ and $\beta_1^i, \dots , \beta_w^i$ be the maximum 
integers, for $j=1,2$, such that $(\xi - \lambda)^{\beta_j^i}$ divides $p_j^i(\xi)$ for 
$j=1,\dots,w$. Then $\beta_j^1 = \beta_j^2$ for $j=1,\dots,w$.
In particular, if
\begin{equation*}
S_1(\xi) = \left[ \begin{array}{cc}
\xi I - P & 0 \\
0 & Q(\xi) \\
\end{array} \right]
\end{equation*}
and 
\begin{equation*}
S_2(\xi) = \left[ \begin{array}{cc}
\xi I - P & 0 \\
0 & I \\
\end{array} \right]
\end{equation*}
and $\det Q(\lambda) \neq 0$, then the partial multiplicities of $S_1(\xi)$ 
and $S_2(\xi)$ corresponding to $\lambda$ are equal.
\end{prop}

Before proving the main result
Theorem \ref{main}, we state and prove another useful result. 
The following lemma relates the partial multiplicities of purely imaginary 
eigenvalues of the Hamiltonian matrix corresponding to the controllable part 
to that of the uncontrollable behavior.

\begin{lemma}\label{evenpartial}
Consider the behavior $\mathfrak{B} \in \mathfrak{L}^\texttt{w}$ with the set 
of uncontrollable modes $\Lambda_{un}$ satisfying $\Lambda_{un} \cap
 -\Lambda_{un} = \emptyset $. Let $\mathfrak{B}$ have an observable 
i/s/o representation as $\frac{d}{dt} x = Ax + Bw_1$, $w_2 = Cx + Dw_1,$ 
induced by $w = (w_1,w_2)$, such that $(I_m+D^TJ_{pq}D) > 0$. Further, let
 $\mathfrak{B}_{cont} = \text{im}\: 
M(\frac{d}{dt})$. Define $\Phi (\zeta, \eta) = M^T(\zeta)\Sigma M(\eta)$ and 
construct the Hamiltonian matrix, $H$ as in \eqref{eq:ham}. Then 
\begin{enumerate}
\item $\sigma(H) = \text{roots}\: (\det \partial \Phi (\xi)) \cup \Lambda_{un} 
\cup (-\Lambda_{un})$.
\item If the controllable part $\mathfrak{B}_{cont}$ is $\Sigma$-dissipative, 
then the partial multiplicities corresponding to the purely imaginary 
eigenvalues of $H$, if any, are all even.
\end{enumerate}
\label{parMul}
\end{lemma}

\noindent \textbf{Proof of Lemma \ref{evenpartial}}:
{\bf Statement 1}: See \cite{PB08}.

\noindent {\bf Statement 2}:
We use the fact that if the controllable part is dissipative, then the partial 
multiplicities corresponding to purely imaginary eigenvalues of $H_c$ are 
even. Without loss of generality the following i/s/o
representation for $\mathfrak{B}$ is assumed

\begin{equation}
\frac{d}{dt} \left[ \begin{array}{c}
	x_c \\
	x_u \\
\end{array}  \right] = 
\left[ \begin{array}{cc}
	\hat{A}_c & \hat{A}_{cp} \\
	0 & \hat{A}_u \\
\end{array}  \right] 
\left[ \begin{array}{c}
	x_c \\
	x_u \\
\end{array}  \right] + 
\left[ \begin{array}{c}
	\hat{B}_c \\
	0 \\
\end{array}  \right] w_1
\label{eq:kal1}
\end{equation}

\begin{equation}
w_2 = 
\left[ \begin{array}{cc}
	\hat{C}_c & \hat{C}_u \\
\end{array}  \right]
\left[ \begin{array}{c}
	x_c \\
	x_u \\
\end{array}  \right] + \hat{D} w_1
\label{eq:kal2}
\end{equation}
with $(\hat{A}_c, \hat{B}_c)$ controllable.

\noindent Then, the Hamiltonian matrix gets the following form


\begin{equation}
H = \left[ \begin{array}{cccc}
	A_c & A_{cp} & B_cB_c^T & 0 \\
	0 & A_u & 0 & 0 \\
	C_c^TC_c & C_c^TC_u & -A_c^T & 0 \\
	C_u^TC_c & C_u^TC_u & -A_{cp}^T & -A_u^T \\
\end{array}  \right]
\label{eq:HamKal}
\end{equation}
where
\begin{equation*}
A_c = \hat{A}_c - 
\hat{B}_c(I_\texttt{m}+\hat{D}^TJ_{pq}\hat{D})^{-1}\hat{D}^TJ_{pq}\hat{C}_c, 
\end{equation*}
\begin{equation*}
A_{cp} = \hat{A}_{cp} - 
\hat{B}_c(I_\texttt{m}+\hat{D}^TJ_{pq}\hat{D})^{-1}\hat{D}^TJ_{pq}\hat{C}_u, 
\end{equation*}
and
\begin{equation*}
B_c = \hat{B}_c(I_\texttt{m}+\hat{D}^TJ_{pq}\hat{D})^{-\frac{1}{2}},
\end{equation*}
\begin{equation*}
C_c = (I_\texttt{pq}+\hat{D}\hat{D}^T)^{-\frac{1}{2}} \hat{C}_c, \: \: \: \: 
C_u = (I_\texttt{pq}+\hat{D}\hat{D}^T)^{-\frac{1}{2}} \hat{C}_u,
\end{equation*}
and $A_u = \hat{A}_u$.

The Hamiltonian matrix for $\mathfrak{B}_{cont}$ is 
\begin{equation*}
H_c = \left[ \begin{array}{cc}
	A_c & B_cB_c^T\\
	C_c^TC_c & -A_c^T \\	
\end{array}  \right]
\end{equation*}

Consider the polynomial matrix $H(\xi) = \xi I_{2n} - H$,
\begin{equation*}
H(\xi) = \left[ \begin{array}{cccc}
	\xi I_{n_c}-A_c & -A_{cp} & -B_cB_c^T & 0 \\
	0 & \xi I_{n_u}-A_u & 0 & 0 \\
	-C_c^TC_c & -C_c^TC_u & \xi I_{n_c}+A_c^T & 0 \\
	-C_u^TC_c & -C_u^TC_u & A_{cp}^T & \xi I_{n_u}+A_u^T \\
\end{array}  \right]
\end{equation*}

Applying the following transformations to $H(\xi)$, we get
$H_1(\xi) = E_1 H(\xi) E_1 = $
\begin{equation*}
\left[ \begin{array}{cccc}
	\xi I_{n_c}-A_c & -B_cB_c^T & -A_{cp} & 0 \\
	-C_c^TC_c & \xi I_{n_c}+A_c^T & -C_c^TC_u & 0 \\
	0 & 0 & \xi I_{n_u}-A_u & 0 \\
	-C_u^TC_c & A_{cp}^T & -C_u^TC_u & \xi I_{n_u}+A_u^T \\
\end{array}  \right]
\end{equation*}
where
\begin{equation*}
E_1 := \left[ \begin{array}{cccc}
	I_{n_c} & 0 & 0 & 0 \\
	0 & 0 & I_{n_c} & 0 \\
	0 & I_{n_u} & 0 & 0 \\
	0 & 0 & 0 & I_{n_u} \\
\end{array}  \right]
\end{equation*}

Now, since $\Lambda_{un} \cap \jj \mathbb{R} = \emptyset$, for $\lambda \in 
\sigma(H_c) \cap \jj  \mathbb{R}$, the matrix blocks $\lambda I_{n_u}-A_u$ and 
$\lambda I_{n_u}+A_u^T$ are invertible.
Pre-multiplying $H_1(\xi)$ by $E_2$ 
and post-multiplying by $E_3$, we get
$H_2(\xi) := E_2 H_1(\xi) E_3 = $
\begin{equation}
\left[ \begin{array}{cccc}
	\xi I_{n_c}-A_c & -B_cB_c^T & 0 & 0 \\
	-C_c^TC_c & \xi I_{n_c}+A_c^T & 0 & 0 \\
	0 & 0 & \xi I_{n_u}-A_u & 0 \\
	0 & 0 & -C_u^TC_u & \xi I_{n_u}+A_u^T \\
\end{array}  \right]
\end{equation}
where
\begin{equation*}
E_2 := 
\left[ \begin{array}{cccc}
	I_{n_c} & 0 & -T_3T_1^{-1} & 0 \\
	0 & I_{n_c} & -T_4T_1^{-1} & 0 \\
	0 & 0 & I_{n_u} & 0 \\
	0 & 0 & 0 & I_{n_u} \\
\end{array}  \right], \: \: \: \:
\end{equation*}
and
\begin{equation*}
E_3 := 
\left[ \begin{array}{cccc}
	I_{n_c} & 0 & 0 & 0 \\
	0 & I_{n_c} & 0 & 0 \\
	0 & 0 & I_{n_u} & 0 \\
	-T_2^{-1}T_4^T & T_2^{-1}T_3^T & 0 & I_{n_u} \\
\end{array}  \right]
\end{equation*}
and $T_1:= \lambda I_{n_u}-A_u$, $T_2:= \lambda I_{n_u}+A_u^T$, $T_3:= -A_{cp}$ 
and $T_4:= -C_c^TC_u$.

Thus
\begin{equation*}
H_2(\xi) = E_2 E_1 H(\xi) E_1 E_2 = \left[ \begin{array}{cc} 
\xi  I_{2n_c} - H_c & 0 \\
0 & Q_u(\xi) \\
\end{array} \right]
\end{equation*}
where 
\begin{equation*}
Q_u = 
\left[ \begin{array}{cc} 
T_1 & 0 \\
-C_u^TC_u & T_2 \\
  \end{array} \right]
\end{equation*}

Thus by using Proposition \ref{evenparmul}, the partial multiplicities of 
purely imaginary eigenvalues of $H$ are even. This completes the proof of 
Lemma \ref{parMul}.

\section{Proof of the main result: Theorem \ref{main}} \label{sec4}

\noindent
\textbf{Proof}:
{\bf If part}: Assume that the controllable part $\mathfrak{B}_{cont}$ is 
$\Sigma$-dissipative and the assumption in the theorem is satisfied. By 
Propositions \ref{LMIunctrb} and \ref{graphneutral}, to prove that the 
behavior $\mathfrak{B}$ is dissipative, it suffices to show the existence of a 
$K \in \mathbb{C} ^{\n \times \n}$ such that the corresponding graph subspace 
is an $n$-dimensional, $M$-invariant, $P$-neutral subspace. To show the 
existence of such a $K$, we use Lemma \ref{invariant} to construct a c-set 
such that the corresponding $n$-dimensional $M$-invariant, $P$-neutral 
subspace is also a graph subspace of $K$. This $K$ is the solution to the ARE 
and storage function for the whole behavior would then be defined as $x^TKx$, 
thus completing the proof.

As the unmixing assumption on uncontrollable modes is
assumed, $\lambda \in \jj \mathbb{R} \cap \sigma (H)$ means
that $\lambda \notin \Lambda_{un}$ and $\lambda \in \sigma(H_c)$. As we
have assumed that the controllable part is $\Sigma$-dissipative, from
Lemma \ref{evenpartial}, the partial multiplicities of real eigenvalues
of $M (:= iH)$ are all even. Using this fact, Lemma \ref{invariant} can be
used to infer that there exists a unique $n$-dimensional $M$-invariant,
$P$-neutral subspace for every c-set.

Now, it remains to show the existence of a c-set such that the corresponding 
$n$-dimensional, $M$-invariant, $P$-neutral subspace is also a graph subspace.

We choose a c-set $\mathcal{C}$ such that $\jj \Lambda_{un} \subseteq 
\mathcal{C}$ and show that the corresponding $n$-dimensional, $M$-invariant, 
$P$-neutral subspace is a graph subspace.
Let $\mathcal{L}$ be the $n$-dimensional, $P$-neutral, $M$-invariant subspace 
of $\mathbb{C}^{2n}$ corresponding to the c-set $\mathcal{C}$ and suppose 
\begin{equation}
\mathcal{L} = \text{Im } \left[  \begin{array}{c}
	X_1 \\
	X_2 \\
\end{array} \right]
\label{eq:L}
\end{equation}
for matrices $X_1$ and $X_2$ $\in \mathbb{C} ^{\n \times \n}$. In order
to prove that $\mathcal{L}$ is a graph subspace it is enough to
prove that $X_1$ is invertible. This is proved using contradiction: we assume
$X_1$ is singular and show that we get a contradiction to the
unmixing assumption on $\Lambda_{un}$. This constitutes the rest of the proof
of the `if part'.  We prove this along the lines of \cite{GLR05}.

$M$-invariance of $\mathcal{L}$ implies that
\begin{equation}
\jj  \left[ \begin{array}{cc}
				\tilde{A} & \tilde{D} \\
				\tilde{C} & -\tilde{A}^\ast \\ 
			\end{array} \right] 
			\left[ \begin{array}{c}
	X_1 \\
	X_2 \\
\end{array}  \right] = \left[ \begin{array}{c}
	X_1 \\
	X_2 \\
\end{array}  \right] T
\nonumber
\end{equation}
for $T \in \mathbb{C}^{\n \times \n}$.
In other words,
\begin{equation}
\jj (\tilde{A}X_1 + \tilde{D}X_2) = X_1T, \qquad \jj (\tilde{C}X_1 - 
  \tilde{A}^\ast X_2) = X_2T.
\label{eq:Hneu1}
\end{equation}
Then as $\mathcal{L}$ is $P$-neutral
\begin{equation}
\left[ \begin{array}{cc}
	X_1^\ast & X_2^\ast \\
\end{array} \right]
\left[ \begin{array}{cc}
	-\tilde{C} & \tilde{A}^\ast \\
	\tilde{A} & \tilde{D} \\
\end{array} \right]
\left[\begin{array}{c}
	X_1 \\
	X_2 \\
\end{array} \right]
= 0, 
\label{eq:neutral1}
\end{equation}
i.e.
$X_2^\ast \tilde{D}X_2 + X_1^\ast \tilde{A}^\ast X_2
+ X_2^\ast \tilde{A}X_1 - X_1^\ast \tilde{C}X_1 = 0.  $


Now suppose $X_1$ is singular. Let $\mathcal{K} = 
\ker X_1$. From \eqref{eq:neutral1}, for every $x \in \mathcal{K}$ we have,
\begin{equation}
x^\ast X_2^\ast \tilde{D}X_2x + x^\ast X_1^\ast \tilde{A}^\ast X_2x + x^\ast 
X_2^\ast \tilde{A}X_1x - x^\ast X_1^\ast \tilde{C}X_1x = 0,
\nonumber
\end{equation}
which implies
\begin{equation}
x^\ast X_2^\ast \tilde{D}X_2x = 0.
\label{eq:Dsemi}
\end{equation}
Since $\tilde{D} \geq 0$, $X_2x \in \ker \tilde{D}$, i.e., $X_2\mathcal{K} 
\subseteq \ker \tilde{D}$.
Now, for every $x \in \mathcal{K}$, from equations \eqref{eq:Hneu1} we have,
\begin{equation}
X_1Tx = \jj \tilde{A}X_1x + \jj \tilde{D}X_2x = 0,
\nonumber
\end{equation}
that is,
\begin{equation}
T\mathcal{K} \subseteq \mathcal{K}.
\label{eq:TKinK}
\end{equation}
This implies $\mathcal{K}$ is $T$-invariant. Hence there exists a non-zero $v$ 
which is an eigenvector of T such that $X_1 v =0$ corresponding to eigenvalue, 
say $\lambda$. We claim that $\lambda$ cannot be a real eigenvalue and prove 
this below.

Post-multiplying the second equation of \eqref{eq:Hneu1} by $v$ we get,

\begin{equation}
\jj \tilde{C}X_1 v - \jj \tilde{A}^\ast X_2 v = X_2T v.
\label{eq:Hneu21}
\end{equation}
\begin{equation}
-\jj  \tilde{A}^\ast X_2 v = \lambda X_2 v.
\label{eq:Hneu22}
\end{equation}
This implies that $X_2 v$ is a left eigenvector of $\tilde{A}$ with eigenvalue 
$-\jj \barlambda$ and from equation \eqref{eq:Dsemi} we have $B^{T}X_2 v = 0$.
 This 
means that $-\jj \barlambda$ is an uncontrollable eigenvalue of $\tilde{A}$ and 
$-\jj \barlambda \in \Lambda_{un}$, i.e., $\barlambda \in \jj \Lambda_{un}$.
 Now, if 
$\lambda$ were a real eigenvalue, then $-\jj \barlambda$ and
$\jj \barlambda$ belong to 
$\Lambda_{un}$ and contradicts the unmixing assumption $\Lambda_{un} \cap 
\Lambda_{un} = \emptyset$.

Since $\sigma(T) \setminus \mathbb{R} = \mathcal{C}$ (by Proposition 
\ref{invariant}), $\lambda \in \mathcal{C}$. We now have $\lambda \in 
\mathcal{C}$ but $\overline{\lambda}$ cannot belong to $\mathcal{C}$. But we 
have $\barlambda \in \jj \Lambda_{un}$ and thus $\barlambda \in \mathcal{C}$
 which is a 
contradiction by definition of c-set. This completes the proof for the if part.

\noindent {\bf Only if part:}
Assume $\mathfrak{B}$ is $\Sigma$-dissipative. Then there exists a storage 
function $Q_{\psi}(w)$ such that 
\begin{equation}
\frac{d}{dt} \, Q_\psi (w) \leqslant Q_{\Sigma} (w), \qquad \forall w \in 
\mathfrak{B}
\end{equation}
Integrating both sides for every $ w \in \mathfrak{B} \cap \mathfrak{D}$,
\begin{equation}
\int_{\mathbb{R}} Q_{\Sigma}(w) dt \geqslant 0
\end{equation}
which implies that $\mathfrak{B}_{cont}$ is $\Sigma$-dissipative.
This completes the proof of Theorem \ref{main}.
\EProof

\noindent
The above proof is constructive in the sense that if a behavior $\B\in\Lw$
satisfies the three conditions:
\begin{itemize}
 \item uncontrollable poles are unmixed, i.e. no two of them add to zero
 \item the controllable part $\Bcont$ is dissipative,
 \item the controllable part $\Bcont$ is `strictly dissipative' at
  infinity, i.e.  $(I_m+D^TJ_{pq}D) > 0$ where $D$ is the feed-through
  term of the transfer function,
\end{itemize}
then we construct a storage function that satisfies the dissipation
inequality for the whole behavior $\B$. Further, the storage function
we construct is equal to $x^T K x$ where $K$ is a solution to
the corresponding Algebraic Riccati Equation. Further, we constructed
the storage function by starting with a state representation
in which the state $x$ is \emph{observable}
from the manifest variable $w$. 
These facts lead to the following important corollary.  

\begin{corollary}
Let $\B\in\Lw$ be an uncontrollable behavior whose uncontrollable poles
$\Lambda_{un}$ are unmixed, i.e.  $\Lambda_{un}\cap -\Lambda_{un} = \emptyset$.
Consider $\Sigma$ partitioned in accordance with the input cardinality of $\B$
as
\[
\Sigma = 
\begin{bmatrix}
        I_{\texttt{m}} & 0 & 0 \\
        0 & I_{\texttt{q}} & 0 \\
        0 & 0 & -I_{\texttt{p}} \\
\end{bmatrix}~.
\]
Suppose $\Bcont$, the controllable part of $\B$, has
an observable image representation, $w= M(\frac{d}{dt}) \ell $,
where $M(\xi)$ is partitioned as 
\begin{equation*}
M(\xi) = 
\begin{bmatrix} W_1(\xi) \\ W_2(\xi) \end{bmatrix}; \: \: \: \: W_1 \in 
\mathbb{R}^{\texttt{m} \times \texttt{m}}\left[ \xi \right], W_2 \in 
\mathbb{R}^{(\texttt{p+q}) \times \texttt{m}}\left[ \xi \right].
\end{equation*} 
Let $G(s) := W_2(s)W_1(s)^{-1}$ and $D := \lim_{s \to \infty} 
G(s)$. Assume $\mathfrak{B}_{cont}$ is such that $(I_m+D^TJ_{pq}D) > 0$.  
Then, the following are equivalent.
\begin{enumerate}
  \item $\Bcont$ is dissipative.
  \item There exists a $\Psi\in \Rwwze$ such that $Q_\Psi(w)$ is a
   storage function, i.e.
    $\der Q_\psi(w)\leqslant w^T \Sigma w$ for all $w \in \B$.
   \item  There exists a matrix $K$ and an observable state variable $x$
   such that $\der x^T K x \leqslant w^T \Sigma w$ for all $w \in \B$.
\end{enumerate}
\end{corollary}
Statement 2 tells that the storage function can be expressed as a quadratic
function of the manifest variables $w$ and their derivatives. Statement
3 says that the storage function is a `state function', i.e. a static
function of the states, and hence storage of energy requires no more memory
of past evolution of trajectories than required for arbitrary concatenation of
any two system trajectories.

\section{Examples} \label{sec5}

In this section we discuss two examples of uncontrollable systems
that are dissipative. The Riccati equations encountered in these
cases are solvable by the methods proposed in this paper; we also 
give solutions to the Riccati equations.

The first example is of an uncontrollable system with 
uncontrollable modes satisfying the unmixing assumption, i.e. no two
of the uncontrollable poles add to zero. However, the Hamiltonian
matrix has eigenvalues on the imaginary axis.

\begin{exa}
Consider the behavior $\B$ whose input/state/output representation is given by 
the following $A,B,C$ and $D$ matrices
\begin{equation} \nonumber
A= \begin{bmatrix} 0 & -0.5 \\
1 & -1.5 \end{bmatrix}, \: \: \: 
B= \begin{bmatrix} -0.5 \\ -0.5 \end{bmatrix}, \: \: \:
C= \begin{bmatrix} 0 & -0.5 \end{bmatrix}, \: \: \: 
D=0.5
\end{equation}
with $\sigma(A) = \{-\frac{1}{2}, -1 \}$. Here $\Lambda _{un} = \{-1 \}$ which 
satisfies the unmixing assumption. An equivalent kernel representation 
of the behavior is given by
\begin{equation}
\begin{bmatrix}
(\xi ^2 + 2\xi + 1) & -(2\xi^2 + 3\xi +1) \\
\end{bmatrix} 
\begin{bmatrix}
w_1 \\
w_2 \\
\end{bmatrix} = 0
\nonumber
\end{equation}
In this case $\Sigma=\mathrm{diag}(1,-1)$ and hence
$\sigma_+ (\Sigma) = \texttt{m}(\mathfrak{B})$, it can 
be checked that the controllable part $\B_{cont} = \ker \begin{bmatrix} (\xi 
+ 1) & -(2\xi + 1) \end{bmatrix}$ is $\Sigma$-dissipative. And 
$(I_m+D^TD) = 3/4 > 0$. Thus from Theorem \ref{main}, $\B$ is 
$\Sigma$-dissipative. 

The following real symmetric matrix induces a storage function that satisfies 
the dissipation inequality
\begin{equation}
K = \frac{1}{6} \begin{bmatrix} 7 & -1 \\ -1 & 1 \\ \end{bmatrix}
\nonumber
\end{equation}

The 2-dimensional, $M$-invariant, $P$-neutral subspace which gives the 
solution is 
\begin{equation}
im \begin{bmatrix}
\jj  & 2 \\
\jj  & 8 \\
\jj  & 1 \\
0 & 1 \\
\end{bmatrix}
\nonumber
\end{equation}
\end{exa}

The next example is an RLC circuit shown in Figure \ref{RLC.exa.fig}.

\begin{exa}
Consider the RLC circuit system whose input $u$ is the current flowing
into the circuit, and output $y$ is the current through the inductor.
A state space representation of the system is found using the
following definition of the states.
The state variables are $x_C$, voltage across the capacitor,
and $x_L$, the current through the inductor.  Assume $R_L=R_C=:R$.

\begin{figure}[!h]
\centerline{\resizebox{!}{35mm}{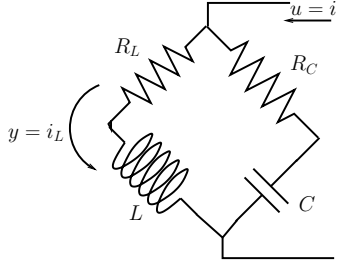}}
  \caption{An RLC circuit system with input $u$ and output $y$}
  \label{RLC.exa.fig}
\end{figure}

\begin{equation*}
A= \begin{bmatrix} 0 & -1/C \\ 1/L & -2R/L \end{bmatrix},
\: \: \: 
B= \begin{bmatrix} 1/C \\ R/L \end{bmatrix}, \: \: \:
C= \begin{bmatrix} 0 & 1 \end{bmatrix}, \: \: \: 
D=R, 
\end{equation*}

The system becomes uncontrollable when $L=R^2C$.  
Let $R=0.5$, $C=1$ and $L=0.25$. The state representation of the system is

\begin{equation}
A= \begin{bmatrix} 
0 & -1 \\
4 & -4\\
\end{bmatrix}, \: \: \: 
B= \begin{bmatrix}
1  \\
2 \\
\end{bmatrix}, \: \: \:
C= 
\begin{bmatrix}
0 & 1
\end{bmatrix}, \: \: \: 
D=0,
\nonumber
\end{equation}

Poles of the system are $-2,-2$ and one of them is uncontrollable.
The controllable part is dissipative and the corresponding
Hamiltonian matrix has eigenvalues on the imaginary axis.
A solution to the ARE is given by the following symmetric matrix
\begin{equation}
K = \begin{bmatrix} 
3 & -0.5 \\
-0.5 & 0.25 \\
\end{bmatrix}
\end{equation} 
which induces the storage function $3x_C^2-x_C x_L + x_L^2/4$.
\end{exa}

\section{Behaviors with static controllable part}

In this section we consider the case when all the states of
the system are uncontrollable, in other words, when the
only controllable part is static/memoryless. In a state space realization,
this means that the matrix $B$ is zero. The case of autonomous behaviors
is clearly a special case: $D$ also is zero and the assumptions in the
lemma are satisfied.

\begin{lemma}
Consider a behavior $\B$ with static controllable part. 
Let the behavior have a state representation,
\begin{equation}
	\frac{d}{dt}x=Ax, \qquad w_2=Cx+Dw_1.
\end{equation}
with the pair $(C,A)$ observable.
Assume $(I_m+D^TJ_{pq}D) > 0$ and $\Lambda _{un} \subset i\R$. Then there
 does not exist a symmetric solution to the corresponding ARE.
\end{lemma}

\noindent
{\bf Proof:}
The Hamiltonian matrix takes the form,
\begin{equation}
H = \begin{bmatrix}
\tilde{A} & 0 \\
\tilde{C} & -\tilde{A}^\ast \\ 
\end{bmatrix} 
\end{equation}

Here $\tilde{A} = A$, $\tilde{C} = C^T(J_{pq}+DD^T)^{-1}C$.

Next we use the following proposition to say that partial multiplicities of
 purely imaginary eigenvalues of the Hamiltonian matrix are even.

\begin{prop}
Consider the matrix
\begin{equation}
N = \begin{bmatrix}
\tilde{A} & \tilde{D} \\
0 & -\tilde{A}^\ast \\ 
\end{bmatrix}.
\end{equation}
Then for every purely imaginary $\lambda _0 \in \sigma(N)$ such that
 $\mathcal{R}_{\lambda_0}(N) \subseteq \mathfrak{C}_{\tilde{A},\tilde{D}}$,
 the partial multiplicities of such $\lambda _0$ are even. In fact, they are
 twice the partial multiplicities of $\lambda _0$ as an eigenvalue of $A$.
\label{part}
\end{prop}

Consider, 
\begin{equation}
H^T = \begin{bmatrix}
\tilde{A}^T & \tilde{C}^T \\
0 & -\tilde{A} \\ 
\end{bmatrix}.
\end{equation}
As the pair $ (\tilde{C},\tilde{A}) $ is observable, the pair
 $(\tilde{A}^T,\tilde{C}^T)$ is controllable. From Proposition \ref{part},
 the partial multiplicities of purely imaginary eigenvalues are twice the
 partial multiplicities of purely imaginary eigenvalues of $A$.

\section{Lossless autonomous behaviors}

In this section we investigate the requirement of unobservable variables
in the definition of the storage function. As has been studied/shown so far,
for controllable dissipative systems (\cite{WT98}), 
the storage function need not depend on unobservable variables.  
 It was later shown in \cite{PB08} that for the case of strict dissipativity
 of uncontrollable systems, observable storage functions are enough under
 unmixing and maximum input cardinality conditions. Theorem \ref{main} shows
 that this is true for a more general scenario, i.e., dissipativity
 (including non-strict dissipativity) for all input cardinality conditions
 under unmixing assumption. On the other hand, when relaxing the 
 unmixing assumption
 elsewhere except on the imaginary axis, under certain conditions solutions
 to the ARI exists though ARE does not have a solution (see \cite{Fai87}).
In this section, we investigate the need for unobservable storage functions for
 uncontrollable systems whose uncontrollable poles lie entirely on the
 imaginary axis. We discuss the case for autonomous behaviors below. 

\begin{lemma}
Consider an autonomous behavior $\B_{aut}$
\begin{equation}
\frac{d}{dt}x = Ax, \: \: \: w=Cx,
\nonumber
\end{equation}
with $\sigma(A) \cap i\mathbb{R} \neq \emptyset $. Let the supply rate be
 $Q_{\Sigma}(w) = -w^Tw$. Then the following is true.

If there exists a storage function $Q_{\psi}(w)$ satisfying the inequality
\begin{equation}
	\frac{d}{dt}Q_{\psi}(w) \leqslant Q_{\Sigma}(w), 
  \qquad \forall \: w \in \B_{aut}
\end{equation}
then any $\lambda \in \sigma(A) \cap i\R$ is C-unobservable.
\end{lemma}

\noindent  \textbf{Proof:}
Suppose if there exists a storage function which is a state function $x^TKx$
 satisfying the dissipation LMI, then the dissipation LMI \eqref{eq:LMInew} is
 equivalent to the Lyapunov inequality
\begin{equation}
KA + A^\ast K + C^\ast C \leqslant 0.
\label{eq:01auto}
\end{equation}
Now for every eigenvector of A corresponding to eigenvalue $\lambda \in
 i\mathbb{R}$, we have
\begin{equation}
x^\ast (KA + A^\ast K + C^\ast C) x \leqslant 0,
\label{eq:03auto}
\end{equation}
which gives
\begin{equation}
\lambda x^\ast Kx + \overline{\lambda} x^\ast Kx + x^\ast C^\ast Cx \leqslant 0
\label{eq:04auto}
\end{equation}
or 
\begin{equation}
x^\ast C^\ast Cx \leqslant 0.
\label{eq:05auto}
\end{equation}
But, as $x^\ast C^\ast Cx \geq 0$, we have $Cx = 0$ for every eigenvector of
 $A$ corresponding to $\lambda \in \sigma(A) \cap i\R$.
This implies that the any $\lambda \in \sigma(A) \cap i\R$ is C-unobservable.

This observation tells that for dissipativity of autonomous systems having
 eigenvalues on the imaginary axis, it is necessary to allow storage
 functions to depend on unobservable variables also.





\section{Orthogonality and uncontrollable behaviors}

In this section we investigate the property of orthogonality of
two behaviors in the absence of controllability. We propose
a definition that is intuitively expected and show that
by relating this definition to lossless uncontrollable behaviors,
we encounter a situation that suggests an exploration
whether dissipativity should
be defined for behaviors for which the input-cardinality condition
is not satisfied.  We first review a result about orthogonality of
controllable behaviors.

\begin{prop} \label{prop:lossless:orthogonal}
Let $\Sigma \in \Rww$ be nonsingular,
and suppose $\B_1, \B_2 \in \Lwcont$. The following are equivalent.
\begin{enumerate}
\item $\int_{\R} w_1 ^T \Sigma w_2 dt = 0$
for all $w_1 \in \B_1 \cap \D$ and for all $w_2 \in \B_2 \cap \D$.
\item $\B_1 \times \B_2$ is lossless with respect
to $\begin{bmatrix} 0 & \Sigma \\ \Sigma ^T & 0 \\ \end{bmatrix}$.
\item There exists a bilinear differential form $L_\Psi$, induced by
$\Psi \in \Rwwze$ such that $\der L_\Psi(w_1,w_2) = w_1^T \Sigma w_2$.
\end{enumerate}
\end{prop}

Statement 1 above is taken as the definition of orthogonality between
two controllable behaviors $\B_1$ and $\B_2$ in \cite{WT98}.
Keeping in line with Definition \ref{defn:dissip} for dissipativity,
we could take Statement 3 above as the definition of orthogonality
for behaviors not necessarily controllable. The drawback of this approach
is elaborated later below in this section. We pursue a different direction as
follows. Notice that if $\B_1$ and $\B_2$ satisfy the integral condition
in Statement 1, then this integral condition is satisfied for every
respective sub-behaviors $\B_1'$ and $\B_2'$ also. Of course, restricting
to compactly supported trajectories in the integration 
implies only controllable parts of respectively $\B_1'$ and $\B_2'$ satisfy
orthogonality. The following definition builds on this property.

\newcommand{\Sperp}{\perp_{\Sigma}}

\begin{definition}
Consider a nonsingular $\Sigma \in \Rww$ and let $\B_1$ and $\B_2 \in\Lw$.
Behaviors $\B_1$ and $\B_2$ are said to be $\Sigma$-orthogonal (and denoted
by $\B_1 \Sperp \B_2$)  if there exist $\B_1^c$ and $\B_2^c \in \Lwcont$
such that 
\begin{itemize}
\item $\B_1 \subseteq \B_1^c$,
\item $\B_2 \subseteq \B_2^c$, and
\item $\int_{\R} w_1 ^T \Sigma w_2 dt = 0$
for all $w_1 \in \B_1^c \cap \D$ and for all $w_2 \in \B_2^c \cap \D$.  
\end{itemize}
\end{definition}

While this definition is not existential in the storage function, it
is existential in $\B_1^c$ and $\B_2^c$, raising questions about
how to check orthogonality. Note that if $\B$ is an uncontrollable
behavior, then any controllable $\B^c \in \Lwcont$ such that
$\B \subseteq \B^c$ satisfies $\m(\B) < \m(\B^c)$, and
$\B \subsetneqq \B^c$.

The question arises as to how much larger a controllable $\B^c$ would have
to be for it to contain $\B$. This problem is addressed in the
following subsection.

\subsection{Smallest controllable superbehavior} \label{subsec:superbehavior}

Due to its significance for determining whether two
uncontrollable behaviors are orthogonal, in this subsection we study
the following problem:

\begin{problem} \label{prob:smallestB}
Let $\B_1 \in \Lw$.  Find $\B_2 \in \Lw$ such that
\begin{enumerate}
    \item $\B_2 \supseteq \B_1 $
    \item $\B_2 \in \Lwcont$ i.e., a controllable behavior
    \item $\B_2$ is a behavior with the smallest input cardinality satisfying
 Properties 1 and 2.
\end{enumerate}
\end{problem}

The following theorem answers this question.

\begin{theorem} \label{thm:superbehavior:controllable}
For $\B \in \Lw$, the following statements are true.
\begin{enumerate}
\item There exists $\B_2 \in \Lwcont$ satisfying
the requirements in Problem \ref{prob:smallestB}.
\item The behavior $\B_2$
is unique if and only if $\B_1$ is controllable, and in that case
$\B_1=\B_2$. 
\item Assume $\B_1$ is uncontrollable.
The input cardinality of $\B_2$, $\m(\B_2)$ satisfies $\m(\B_2)> \m(\B_1)$. 
More precisely, $\m(\B_2)=\m(\B_1)+\kt$ where 
$\kt:=\underset{\lambda\in\C}{\max}
  \rank(R_1(\lambda))-\underset{\lambda\in\C}{\min} 
\rank(R_1(\lambda))$ where $R_1$ is a kernel representation of $\B_1$. 
\end{enumerate} 
\end{theorem}

\begin{proof1}
{\bf (1):} This is shown by constructing such a $\B_2$.
Let $\ker R_1(\xi)$ and $\ker R_2(\xi)$ be the kernel representations 
of $\B_1$ and $\B_2$ respectively where $R_1 \in \mathbb{R}^{p_1 \times w}$, 
$R_2 \in \mathbb{R}^{p_2 \times w}$ and there exists an 
$F \in \mathbb{R}^{p_2 \times p_1}$ such that 
$FR_1 = R_2$. Let $\Lambda_{un}$ be the set of uncontrollable modes of $\B_1$.
Without loss of generality, we assume $R_1$ to be 
\begin{equation}
R_1 = \left[ \begin{array}{cc}
	S & 0 
\end{array} \right]
\end{equation}
where $S$ is of the form $\left[ \begin{array}{cc}
	I & 0 \\
	0 & D 
\end{array} \right]$ such that $I \in \mathbb{R}^{(p-\kt) \times (p-\kt)}$ is
 identity matrix and $D$ $\in \mathbb{R}^{\kt \times \kt}[\xi]$ is a diagonal
 matrix with $d_1, d_2 \dots, d_\kt$ along its diagonal and satisfying the
 divisibility property: $d_1 | d_2$, $d_2 | d_3,$,  $\dots$ and $d_{\kt-1} | d_\kt$,
 with degree of $d_1$ at least one, and with
\[
\kt=\underset{\lambda\in\R}{\max}  \rank(R_1(\lambda)) - 
   \underset{\lambda\in\R}{\min} 
\rank(R_1(\lambda)). 
\]

\noindent Partitioning $F = \left[ \begin{array}{cc}
	F_1 & F_2
\end{array} \right]$ conforming to the row partition of $R_1$, we have
\begin{equation}
R_2 = \left[ \begin{array}{cc}
	F_1 & F_2
\end{array} \right]
\left[ \begin{array}{ccc}
	I & 0 & 0 \\
	0 & D & 0
\end{array} \right]
\end{equation}
This simplifies to
\begin{equation}
R_2 = \left[ \begin{array}{ccc}
	F_1 & F_2D & 0 \\	
\end{array} \right]
\end{equation}
Let $\B_2$ be the behaviour defined by the kernel representation of $R_2$.
 For $\B_2$ to be controllable, $R_2(\lambda)$ needs to have full row rank
 for all $\lambda \in \mathbb{C}$. As $F_2P$ loses rank for 
$\lambda \in \Lambda_{un}$, $F_1$ should have full row rank for every
 $\lambda \in \C$ so that $\B_2$ is controllable. This proves the existence
 of $\B_2$ satisfying properties 1 and 2 of \ref{prob:smallestB}.
In order to satisfy property 3 in Problem \ref{prob:smallestB}, 
i.e. $\texttt{m}(\B_2)$ has to be the least, we have to choose a 
unimodular $F_1$ and free $F_2$ such that $\B_2$ satisfies the three 
conditions. $F_2$ can be freely chosen because the choice of $F_2$ does 
not affect the input cardinality of $\B_2$.

\noindent
{\bf (2):}
Let $\B_1$ be controllable. Then, in the above, $S=I$ and $D$ does not exist.
 This means that $F_2$ does not exist. Thus the kernel representation matrix
 of $\B_2$ would be given by $[F_1 \hspace{5pt} 0]$ where $F_1$ is
 unimodular. Therefore, $\B_2$ is unique.
Let $\B_1$ be uncontrollable. Then $D$ exists and as shown above, free $F_2$
 can be chosen. This makes the behaviour $\B_2$ non-unique. 

\noindent
{\bf (3):} The input cardinality of $\B_2$ is given by
\begin{equation}
\begin{split}
\m(\B_2) & = w - p_2 \\
& = w - (p_1 - \kt) \\
& = \m(\B_1) + \kt
\end{split}
\end{equation}
Here, $\kt=\underset{\lambda\in\R}{\max}  \rank(R_1(\lambda)) - 
 \underset{\lambda\in\R}{\min} \rank(R_1(\lambda))$ which is determined by
 the size of the $D$ matrix. 
\end{proof1}

\subsection{Superbehaviors and orthogonality}

We saw above in this section that orthogonality of two uncontrollable behaviors
is defined by requiring these uncontrollable behaviors to be sub-behaviors
of two orthogonal controllable behaviors. Using the result on existence
of super-behaviors that are controllable, and their non-uniqueness
even if the are the smallest controllable superbehavior, we formulate
the question of whether two uncontrollable behaviors are orthogonal
as a question of finding a pair of smallest controllable superbehaviors that
are mutually orthogonal. The requirement of them being smallest is motivated
by the fact that orthogonality of two controllable behaviors imposes
an upper bound on their input cardinalities: this is reviewed
below.   For a behavior $\B \in \Lw$ and a nonsingular matrix $\Sigma$,
the set $\Sigma\B$ is defined
as follows
\[
\Sigma \B := \{ w \in \CinfRRw \mid \mbox{ there exists } v \in \B
\mbox{ such that } w=\Sigma v\}.
\]
It is straightforward that $\Sigma \B$ is also a behavior,
its controllability is equivalent to that of $\B$, and the input cardinalities
are equal.

\begin{prop}
Let $\B_1$ and $\B_2\in\Lwcont$ and suppose $\Sigma \in \Rww$ is nonsingular.
Then, the following are true.
\begin{enumerate}
\item $\B_1 \Sperp \B_2$ $\Leftrightarrow$ $\B_1 \perp (\Sigma \B_2)$.
\item $\B_1 \perp \B_2 $ $\Rightarrow$ $\m(\B_1)+\m(\B_2) \leqslant \wt$.
\item $\B_1 \Sperp \B_2$ $\Rightarrow$ $\m(\B_1)+\m(\B_2) \leqslant \wt$.
\end{enumerate}
\end{prop}

Due to the above inequality constraint on the input-cardinalities
of orthogonal controllable behaviors, the uncontrollable
behaviors too have a necessary condition to satisfy for
mutual orthogonality.

\begin{lemma} \label{lem:uncontrol:ortho}
Suppose $\B_1$ and $\B_2 \in \Lw$ with at least one of them uncontrollable 
and let $\Sigma \in \Rww$ be nonsingular.
Assume $\B_1 \Sperp \B_2$. Then $\m(\B_1)+\m(\B_2) < \wt$.  
\end{lemma}

\begin{example}
Consider the pair of `seemingly' orthogonal behaviors
studied in \cite[page 360]{Wil04}. Define $\B_1$ and $\B_2 \in \Lw$ by 
\[
\der x = A x, \quad w_1 = Cx, \qquad w_2 \mbox{ free, i.e. $B_2 = \CinfRRw$ }
\]
and the supply rate $w_1^T w_2$. Thus $\B_1$ is autonomous. Consider
the following non-observable latent variable representation for $\B_2$:
$\der z = -A^T z + C^T w_2$. It can be checked that the
`storage function' $x^T z$ satisfies $\der x^T z = w_1 ^T w_2$. Existence
of such a storage function is, in fact, reasonable for a (different) definition
of orthogonality of the two behaviors $\B_1$ and $\B_2$.  In fact, any autonomous
behavior $\B_1\in\Lw$ is then `orthogonal' to $\B_2=\CinfRRw$!
However, the `embeddability' definition we have used above rules out
this example for an orthogonal pair of behaviors since
the necessary
condition of the above lemma is not satisfied. In other words, there
doesn't exist a controllable behavior $\B_1^c$ such that $\B_1^c$ contains
$\B_1$ and $\B_1^c \Sperp \B_2$. Thus $\B_1$ and $\B_2$ are not $I$-orthogonal.
\end{example}

\section{Dissipative sub-behaviors/superbehaviors}

In this section we look into the input cardinality condition for dissipative
behaviors. Recall that a behavior $\B \in \Lw$ which is dissipative
with respect to the supply rate $\Sigma$ (constant, symmetric, nonsingular
matrix) satisfies the condition $\m(\B) \leqslant \sigma_+(\Sigma)$.
We now look into the possibility of embedding a behavior $\B$ 
in a controllable superbehavior that is $\Sigma$-dissipative, and
into the drawback of using this as the definition of dissipativity, along
the lines of orthogonality defined in the previous section.

\begin{problem} \label{prob:plusminus}
Given a nonsingular, symmetric and indefinite $\Sigma \in \Rww$, find
conditions for existence of a behavior $\B\in\Lw$ such that
\begin{itemize}
\item there exist $\B_+$ and $\B_- \in \Lwcont$ with $\B=\B_+\cap\B_-$.
\item $\B_+$ is strictly $\Sigma$ dissipative
\item $\B_-$ is strictly $-\Sigma$ dissipative.
\end{itemize}
\end{problem}
The significance of the above problem is that if a nonzero
behavior $\B$ satisfying
above conditions exists, then clearly such a behavior would be
both strictly $\Sigma$ and strictly $-\Sigma$ dissipative, raising concerns
about whether embeddability in a dissipative controllable superbehavior
is a reasonable definition of dissipativity (when dealing with
uncontrollable behaviors). The following theorem states that
nonzero autonomous behaviors can indeed exist satisfying above condition.

\begin{theorem} \label{thm:strict:embeddable:autonomou}
Let $\Sigma \in \Rww$ be nonsingular, symmetric and indefinite.
Then there exists a nonzero $\B\in\Lw$ such that requirements in
Problem \ref{prob:plusminus}
are satisfied. Further, any such $\B$ satisfies
$\m(\B)=0$, i.e. $\B$ is autonomous. 
\end{theorem}

\begin{proof1}
We first show that given any nonsingular, symmetric and indefinite
$\Sigma \in \Rww$, a behavior $\B$ satisfying above properties exists. Without
loss of generality, let 
\[
\Sigma = \begin{bmatrix} I_+ & 0 \\ 0 & -I_- \end{bmatrix}
\]
with sizes of the identity matrices $I_+$ and 
$I_-$ equal to $\sigma_+$ and $\sigma_-$.
Due to indefiniteness of $\Sigma$, 
$\sigma_+,\sigma_- \geqslant 1$ and, due to the nonsingularity, 
they sum up to $w$.  Choose any\footnote{Such matrices are plenty 
due to existence of sufficiently many controllable
strictly dissipative behaviors for every supply rate $\Sigma$ 
satisfying the stated conditions.} nonzero polynomial matrix 
$M_+\in\R^{w \times m_+}[\xi]$ with 
$M_+^T(-\jomega) \Sigma M_+(\jomega) \geqslant \epsilon_+ I _w$ for all
 $\omega \in \R$
for some $\epsilon_+ > 0$. 
 Similarly, 
choose $M_- \in \R^{w \times m_-}[\xi]$ such that
$M_-^T(-\jomega) \Sigma M_-(\jomega) \leqslant -\epsilon _- I_w $ for 
all $\omega \in \R$ for some $\epsilon_- > 0$. 
Now define $\B_+$ and $\B_-$ by image representations
$w=M_+(\der)\ell$ and $w=M_-(\der)\ell$ respectively. Define 
$\B:=\B_+\cap \B_-$, thus proving\footnote{It is 
not difficult to show that if $M_+$ and $M_-$ were nonconstant 
polynomial matrices, then $\B$ is not the zero behavior. An example
following the proof makes this easier to see.}
existence of $\B$ as stated in the theorem.

We now show that $\B$ is autonomous.
Let $R_+$ and $R_-$ be minimal kernel representation matrices of 
$\B_+$ and $\B_-$ respectively. Then $\B=\B_+\cap\B_-$ is described by 
the kernel representation matrix 
$R\in \R^{\bullet \times w}[\xi]$ with 
$R :=\left[ \begin{array}{ c } R_+ \\ R_- \end{array} \right]$ .
Clearly, $\rank(R)\leqslant w$. 
Suppose $\rank(R)< \wt $. 
Then, there exists $p\in \R^{\wt}[s]$ and $p \neq 0$ such 
that $R_+p=0$ and $R_-p=0$. This implies $\Image(p(\frac{d}{dt})) \in \B_+$ 
and $\Image(p(\frac{d}{dt})) \in \B_-$. Taking $w=p(\frac{d}{dt})\ell$ with 
$\ell \in \D(\R,\R)$ and $\ell \neq 0$, it follows that 
$w \in \D(\R^{\wt},\R)$. Further, $p\neq0$, hence $w\neq0$ 
because $\ell$ is non-zero and of compact support. 
Further, we have $\epsilon_+,  \epsilon_- > 0$ such that
\[
\int_{-\infty}^{\infty} w^T \Sigma wdt \geqslant \epsilon_+ \Vert 
w \Vert_{L_2}^{2}
\qquad
\mbox{and}
\qquad
\int_{-\infty}^{\infty} w^T \Sigma wdt \leqslant -\epsilon_- \Vert w \Vert_{L_2}^{2}
\]
Both the above conditions cannot be satisfied simultaneously for $w\neq0$. 
Thus, $\rank(R) < \wt$ gives a contradiction. This proves $\rank(R)=\wt$ and 
hence autonomy of $\B$.  
\end{proof1}

We illustrate the above theorem using an example.  
\begin{example} Let $\Sigma = \diag(1,-1)$. Define $\B_+$ and $\B_-$ by
image representations $w=M_+(\der)\ell$ and
$w=M_-(\der)\ell$ respectively, with
\[
M_+(\xi) = \begin{bmatrix} \xi+4 \\ 3 \end{bmatrix} \quad \mbox{ and }
M_-(\xi) = \begin{bmatrix} 2 \\ \xi+5 \end{bmatrix}.
\]
Strict dissipativities is easily verified. Calculating the kernel 
representations, we get a kernel representation for $\B:=\B_+\cap \B_-$
as $R(\der)w = 0$
\[
R(\xi) = \begin{bmatrix} -3 & \xi + 4 \\ \xi+5 & -2 \end{bmatrix}.
\]
Clearly, $R$ is nonsingular and hence $\B$ is autonomous.
\end{example}

For non-strict dissipativity case, we have the following problem and theorem.

\begin{problem} \label{prob:plusminus2}
Given a nonsingular, symmetric and indefinite $\Sigma \in \Rww$, find
conditions for existence of a behavior $\B\in\Lw$ such that
\begin{itemize}
\item there exist $\B_+$ and $\B_- \in \Lwcont$ such that $\B=\B_+\cap\B_-$.
\item $\B_+$ is $\Sigma$ dissipative
\item $\B_-$ is -$\Sigma$ dissipative.
\end{itemize}
\end{problem}

\begin{theorem}
Let $\Sigma \in \Rww$ be nonsingular, symmetric and indefinite.
Then there exists $\B\in\Lw$ such that requirements in
Problem \ref{prob:plusminus2}
are satisfied. Any such $\B$ satisfies
$\m(\B) \leqslant \min(\sigma_+(\Sigma),\sigma_-(\Sigma))$. In 
case $\B$ is uncontrollable, $\m(\B) < \min(\sigma_+(\Sigma),\sigma_-(\Sigma))$\\
If $m(\B) \geqslant 1$, then neither $\B_+$ nor $\B_-$ can be strictly dissipative.
\end{theorem}

\begin{proof1}
The proof proceeds in the same way as the proof for the previous theorem, 
except for the strictness of the dissipativities.  Construct $\B_+$
and $\B_-$ as in the previous proof, but with $\epsilon_+$ and 
$\epsilon_-$ equal to zero. We have
\[
\int_{-\infty}^{\infty} w^T \Sigma wdt \geqslant 0 \mbox{ for all } \quad w \in
 \B_+ \cap \D 
\qquad \mbox{ and } \qquad
\int_{-\infty}^{\infty} w^T \Sigma w dt \leqslant 0 \mbox{ for all } \quad w \in
 \B_- \cap \D 
\]
The above two equations imply $ \int_{-\infty}^{\infty} w^T \Sigma wdt = 0$ 
for all $w \in \B \cap \D $.
Since $\B = \B_+ \cap \B_-$, the behavior $\B$ is dissipative with 
respect to both $\Sigma$ and $-\Sigma$. Dissipativity with respect 
to $\Sigma$ implies $m(\B) \leqslant \sigma_+(\Sigma)$. 
Similarly, dissipativity with respect to $-\Sigma$ 
implies $m(\B) \leqslant \sigma_-(\Sigma)$.
This implies
\begin{equation} \label{dsub}
m(\B) \leqslant \min(\sigma_+,\sigma_-)
\end{equation}
If $\B$ is uncontrollable, then from 
Theorem \ref{thm:superbehavior:controllable}, then the two 
inequalities leading to the inequality \eqref{dsub} are both strict. 
Hence, the input cardinality of $\B$ has to be strictly less than that of
 $\B_+$ and $\B_-$. This implies
\begin{equation} 
m(\B) < \min(\sigma_+,\sigma_-)
\end{equation}

If $m(\B) \geqslant 1$, then from Theorem \ref{thm:strict:embeddable:autonomou},
$\B_+$ and $\B_-$ cannot be strictly dissipative with respect 
to $\Sigma$ and $-\Sigma$ respectively. This completes the proof.
\end{proof1}

As one of the consequences of the above theorem,
if the input cardinality condition is satisfied
for an uncontrollable behavior, i.e. $\m(\B)=\sigma_+(\S)$ or 
$\m(\B)=\sigma_-(\S)$, then such a behavior
cannot be embedded into both a $\S$-dissipative controllable behavior and
a $-\S$-dissipative controllable behavior.  However,
an observable storage function for such a situation exists when
the controllable part is strictly dissipative at $\infty$ and when
the uncontrollable modes satisfy the unmixing condition 
(see \ref{main} and \cite{KariBelur:ifac11, PB08} for this
situation in presence of more assumptions),
A situation when 
$\m(\B)=\sigma_+(\S)$ is very familiar: we deal with RLC circuits in the
next section.

6. We use this method of defining orthogonality of two behaviors
to explore further the definition of dissipativity of a behavior.
Here we bring out a fundamental significance of the so-called input-cardinality
condition: the condition that the number of inputs to the system is
equal to the positive signature of the matrix that induces the power supply.
We show that when this condition is not satisfied, then a behavior
could be both supplying and absorbing net power, and is still not lossless.

%

\section{RLC realizability}

In this brief section we revisit a classical result: a rational
transfer function matrix
being positive real is a necessary and sufficient condition for that
transfer matrix to be realizable using only resistors, capacitors and
inductors (see \cite{Bru:31} and also \cite{BotDuf:49} for the case with
transformers).  Note that the transfer matrix captures only
the controllable part of the behavior and positive realness of
the transfer matrix is nothing but dissipativity with respect to
the supply rate $v^T i$. This is made precise below.

Consider an $n$-port electrical network (with each port having two
terminals) and the variable $w=(v,i)$, where $v$ is the vector of voltages 
across
the $n$-ports and $i$ is the vector of currents through these ports, with
the convention that $v^T i$ is the power flowing {\em into} the network.
Behaviors that are dissipative with respect to the  supply rate $v^T i$
are also called\footnote{It is also common to require the storage
function in the context of this dissipativity to satisfy non-negativity.
The sign-definiteness of the storage function is not the focus of
this paper, hence we ignore this aspect.} `passive'.
Given a controllable behavior $\B$ whose transfer matrix $G$ with
respect to a specific input/output partition, say current $i$ is the
input and voltage $v$ is the output, is positive real, one
can check that this behavior is passive. $G$ in this case
is the impedance matrix and is square, i.e. the number of inputs
is equal to the number of outputs. One can introduce
additional laws that the variables need to satisfy, thus resulting
in a sub-behavior $\Bsub\subsetneqq \B$ 
which has a lesser number of inputs;  consider for example these additional
laws as putting certain currents equal to zero: due to opening of
certain ports.  
However, the transfer matrix for
$\Bsub$ with respect to the input/output partition: input as currents through the non-open
ports and output as the voltages across {\em all} the ports,
 is clearly not square,
and in fact, tall, i.e. has strictly more rows than columns.
Let $\Gsub$ denote this transfer function. Since $\Bsub \subset \B$,
the behavior $\Bsub$ is also passive. Of course, $\Bsub$ need not
be controllable, even if $\B$ is assumed to be controllable.
As an extreme case, suppose all the currents are equal to zero, and
we obtain an autonomous $\Bsub$. While RLC realization of such
transfer matrices which are tall, and further of autonomous behaviors
obtained by, for example, opening all ports, has received hardly
any attention, we remark here one very well-studied sub-behavior of
every passive behavior: the zero behavior. 

Consider the single port network with $v=0$ and $i=0$. This port which
is called a nullator (also nullor, in some literature) behaves
as both the open circuit and shorted circuit: see \cite[page 75]{Bel68} 
and \cite{Car64}. The significant fact about a nullator is that
a nullator cannot be realized using only passive elements, and moreover
any realization (necessarily active) leads to the realization of both a nullator
and its companion, the norator; a norator is a two-terminal port
that allows both the voltage across it and current through it to
be arbitrary.


%
%

\section{Concluding remarks}

We briefly review the main results in this paper. We first used
the existence of an observable storage function as the definition
of a system's dissipativity and proved that 
the dissipativities of a behavior and its controllable part are
equivalent assuming the uncontrollable poles are unmixed and
the dissipativity at infinity frequency is strict. This result's proof
involved new results in the solvability of ARE and used indefinite linear
algebra results.

We showed that for lossless autonomous dissipative systems, the storage function
cannot be observable, thus motivating the need for unobservable storage
functions. We then studied orthogonality/lossless behaviors in the
context of using the definition of existence of a controllable
dissipative superbehavior as a definition of dissipativity.
In addition to results about the smallest controllable superbehavior,
we showed necessary conditions on the number of inputs 
for embeddability of lossless/orthogonal behaviors in larger controllable
such behaviors.

In the context of embeddability as a definition of
dissipativity, we showed that one can always
find behaviors that can be embedded in both a strictly dissipative
behavior and a strictly `anti-dissipative' behavior, thus raising
a question on the embeddability definition.  We related this
question to the well-known result that
the nullator one-port circuit is not realizable using only RLC elements.

%
%
%


\begin{thebibliography}{20}

\bibitem{Ath71} M. Athans, The role and use of the stochastic 
Linear-Quadratic-Gaussian problem in control systems design, \emph{IEEE 
Transactions on Automatic Control}, vol. 16, pages 529-552, 1971.

\bibitem{Bel68} V. Belevitch, {\em Classical Network Theory},
Oakland, CA: Holden-Day, 1968.


\bibitem{BLW91} S. Bittanti, A. J. Laub and J. C. Willems (Eds.), \emph{The 
Riccati Equation}, Springer-Verlag, 1991.

\bibitem{BotDuf:49} R. Bott and R. J. Duffin, 
Impedance synthesis without transformers, {\em Journal of  Applied Physics}, 
vol. 20, page 816, 1949.

\bibitem{BoyGhaFerBal:94} S. Boyd, L. E. Ghaoui, E. Feron
 and V. Balakrishnan, {\em Linear Matrix Inequalities in
    System and Control Theory}, SIAM, Philadelphia, 1994.


\bibitem{Bru:31} O. Brune,  Synthesis of a finite two-terminal 
network whose driving-point impedance is a prescribed function of 
frequency, {\em Journal of Mathematics and Physics,} 
vol. 10, pages 191–236, 1931.


\bibitem{CWB03} M.K. \Camlibel, J. C. Willems and M. N. Belur, On the 
dissipativity of uncontrollable systems, in \emph{Proceedings of the 42nd IEEE 
Conference on Decision and Control}, Hawaii, USA, December 2003.

\bibitem{Car64} H.J. Carlin, ``Singular network elements," {\em IEEE 
Transactions on Circuit Theory}, vol. 11, no. 1, pages 67-72, 1964.

\bibitem{Fai87} L. E. Faibusovich, {Matrix Riccati inequality: existence 
of solutions}, \emph{Systems \& Control Letters}, vol. 9, pages 59-64,
 January 1987.

\bibitem{GLR05} I. Gohberg, P Lancaster, and L. Rodman, \emph{Indefinite Linear 
Algebra and Applications}, Birkhauser, Basel, 2005.

\bibitem{Mac63} A. G. J. MacFarlane, {An eigenvector solution of the 
optimal linear regulator problem},  \emph{Journal of Electronics and Control},
vol. 14, no. 6, pages 643-654, June 1963.

\bibitem{PB08} D. Pal and M. N. Belur. {Dissipativity of uncontrollable 
systems, storage functions, and Lyapunov functions}, \emph{SIAM Journal on
Control and Optimization}, vol. 47, no. 6, pages 2930-2966, 2008.

\bibitem{PilSha:98} H.K. Pillai and S. Shankar, 
A behavioural approach to control of distributed systems, {\em SIAM
Journal on Control and Optimization}, vol. 37, pages 388–408,  1998. 

\bibitem{Pot66} J. E. Potter, {Matrix quadratic solutions}, \emph{SIAM Journal 
of Applied Mathematics}, vol. 14, no. 3, pages 496-501, 1966.

\bibitem{RW97} P. Rapisarda and J. C. Willems, {State maps for linear 
systems}, \emph{SIAM Journal on Control and Optimization}, vol. 35, no. 3, 
pages 1053-1091, 1997.

\bibitem{KariBelur:ifac11} S. Karikalan and M.N. Belur, 
Uncontrollable dissipative dynamical systems, in {\em Proceedings of 
the IFAC World Congress}, Milan, Italy, 2011.

\bibitem{Sch:90} C.W. Scherer, {\em The Riccati Inequality and 
State-space $\hinf$-optimal Control}, Ph.D. thesis,
University of W\"urzburg, Germany, 1990.

\bibitem{TW97} H. L. Trentelman and J. C. Willems, {Every storage 
function is a state function}, \emph{Systems \& Control Letters}, vol. 32, 
pages 249-259, 1997

\bibitem{Wil:72} J. C. Willems, Dissipative dynamical systems - 
  Part I: General theory, Part II: Linear systems with quadratic supply rates, 
  {\em Archive for Rational Mechanics and Analysis}, vol. 45, 
  pages 321–351, 352–393, 1972,




\bibitem{WT98} J. C. Willems and H. L. Trentelman, {On quadratic 
differential forms}, \emph{SIAM Journal on Control and Optimization,}
vol. 36, no. 5, pages 1703-1749, 1998.

\bibitem{WT02} J.C. Willems and H.L. Trentelman,
Synthesis of dissipative systems using quadratic differential
forms, \emph{IEEE Transactions on Automatic Control}, vol. 47, pages 53-69, 2002.

\bibitem{Wil04} J. C. Willems, Hidden variables in dissipative systems,
in \emph{Proceedings of the 43rd IEEE 
Conference on Decision and Control}, Bahamas, pages 358-363, 2004.


\bibitem{Wil07-1} J. C. Willems, {Dissipative dynamical systems},
\emph{European Journal of Control}, vol. 33, pages 134-151, 2007.

\bibitem{Wil07-2} J.C. Willems,
The behavioral approach to open and interconnected systems, \emph{Control
Systems Magazine}, vol. 27, pages 46-99, 2007.

\bibitem{Won79} W. M. Wonham, \emph{Linear Multivariable Control: a Geometric 
Approach}, Springer-Verlag, New York, 1979.

\end{thebibliography}
\end{document}